\documentclass[a4paper,12pt]{amsart}
\usepackage[utf8]{inputenc}
\usepackage[T1]{fontenc}
\usepackage[english]{babel}
\usepackage{color}
\usepackage{amsmath}
\usepackage{amssymb}
\usepackage{amsthm}
\usepackage{graphicx}
\usepackage[colorlinks=true, linkcolor=blue]{hyperref}
\usepackage{geometry}
\usepackage{enumitem}
\usepackage{lmodern}
\usepackage{microtype}
\usepackage{subfiles}
\usepackage{soul}
\usepackage[capitalise]{cleveref}
\usepackage{etoolbox}

\makeatletter
\patchcmd{\@setauthors}
  {\MakeUppercase{\authors}}
  {\authors}
  {}{}
\makeatother
\usepackage[colorinlistoftodos,prependcaption,textsize=tiny]{todonotes}
\newcommand {\R} {\mathbb{R}}
\newcommand {\Z} {\mathbb{Z}}

\newcommand {\p} {\partial}

\newcommand{\eps}{\varepsilon}

\DeclareMathOperator{\dive }{div}

\DeclareMathOperator{\supp }{supp}

\newtheorem{theorem}{Theorem}

\newtheorem{thm}[theorem]{Theorem}
\newtheorem{lemma}[theorem]{Lemma}
\newtheorem{prop}[theorem]{Proposition}
\newtheorem{cor}[theorem]{Corollary}

\theoremstyle{definition}
\newtheorem{definition}[theorem]{Definition}
\newtheorem{remark}[theorem]{Remark}

\newtheorem{ass}[theorem]{Assumption}

\date{}

\begin{document}
\title{Dispersion of Magnetic-Coriolis waves}
\author{Niklas Knobel}\thanks{
Department of Mathematics, Imperial College London,
London SW7 2AZ, United Kingdom
nknobel@ic.ac.uk},
\author{Zhuoya YAO}\thanks{
Institut Camille Jordan, Universit\'e Claude Bernard Lyon 1,
Villeurbanne, 69100, France
yao@math.univ-lyon1.fr}
\begin{abstract}
We consider the $3$-D rotating magnetohydrodynamic (MHD) equations. The combination of rotation and background magnetic field leads to the occurrence of dispersive magneto-Coriolis waves in the linearized equations. We establish a $t^{-1}$ decay estimate for waves with dispersion relation $\sigma (\xi)=\xi_3 h(\vert\xi\vert)$ under suitable assumptions on $h$, and show that this includes magneto-Coriolis waves. Then we apply these dispersive estimates to the nonlinear rotating MHD equations and obtain an improved lower bound on the lifespan of solutions.
\end{abstract}

\maketitle
\setcounter{tocdepth}{1}
\tableofcontents
\section{Introduction}

Rotating stellar objects, such as stars and the liquid metal core of planets, are permeated by both rapid rotation and strong magnetic fields. The interaction between the Coriolis and Lorentz forces organizes the dynamics into magneto-Coriolis waves, which propagate along the rotation axis and field lines. While the dispersion relations of these waves have been known in physics since the mid‑20th century, a mathematical study is lacking.  The purpose of this paper is to develop
such a dispersive framework for the rotating MHD equations. The main analytical ingredient is a general dispersive estimate
for phases of the form $\sigma(\xi)=\xi_3 h(|\xi|).$
Under appropriate nondegeneracy assumptions on $h$, we prove that
the linear evolution exhibits $t^{-1}$ decay in three
dimensions. Then we verify that the dispersion relations of the
magneto-Coriolis waves satisfy these assumptions. Finally, we
combine the resulting dispersive estimates with nonlinear energy
estimates to obtain an improved lower bound for the lifespan of
small-data solutions.

We study this phenomenon through the ideal rotating MHD equations around a constant magnetic field

\begin{align}
\begin{split}
\p_t u +(u\cdot \nabla )u  + 2\Omega e_3 \wedge u  + \nabla \pi &= ((b+B_0)\cdot\nabla) b,\\
\p_t b +(u \cdot\nabla) b &=( (b+B_0)\cdot\nabla) u,\label{eq:nl}\\
\dive(u)&=\dive(b)=0,\\
(u(0,x), b(0,x))&=(u_0(x), b_0(x)).
\end{split}
\end{align}

This equation models a conducting fluid rotating about the $e_3$-axis with angular velocity $\Omega\in \R$ in a homogeneous magnetic field  $B_0\in \R^3$ \cite{davidson_2016,lehnert1954magnetohydrodynamic, lehnert1955magnetohydrodynamic}. The velocity field is denoted by $u:\ \R_+\times \R^3 \to \R^3 $ and the magnetic field by $b:\ \R_+\times \R^3 \to \R^3 $. When $\Omega=0$, this reduces to the usual ideal MHD equations.

The mathematical study of MHD began with Alfvén’s pioneering work \cite{alfven1942existence}, which identified Alfvén waves in a conducting fluid permeated by a uniform magnetic field. Alfvén waves are transverse plasma waves that propagate along a background magnetic field. When rotation is included, the Coriolis force splits the Alfvén-wave dynamics into two magneto-Coriolis branches. In particular, while Alfvén waves are not dispersive, the magneto-Coriolis waves are dispersive.  To understand the behavior of equation \eqref{eq:nl}, we consider the linearized equations 
\begin{align}\begin{split}
\p_t u + 2\Omega e_3 \wedge u  + \nabla \pi &=(B_0 \cdot\nabla )b,\\
\p_t b&=(B_0\cdot  \nabla) u,\\
\dive(u) &=\dive(b)=0.\label{eq:lin}
\end{split}\end{align}

The linearized system admits the explicit representation:

\begin{prop}
For the initial data
$u_{0},b_{0}\in H^N(\mathbb R^3)$,
there exist $\alpha_j,\beta_j$,
such that  the global in time solution $(u,b)$ of \eqref{eq:lin} is given by 
\begin{align*}
    u(t,x)&= \sum_{j=1}^4 \int \alpha_j(\xi) \exp(i(\xi \cdot x - \sigma_j(\xi) t ))d\xi,\\
    b(t,x)&= \sum_{j=1}^4 \int \beta_j(\xi) \exp(i(\xi \cdot x - \sigma_j(\xi) t ))d\xi,
\end{align*}

with the dispersion relation
\begin{align}\label{sigmaj}
    \sigma_{1}(\xi)&=\Omega\frac {\xi_3}{\vert \xi \vert } + \sqrt{ \Omega^2\frac {\xi_3^2}{\vert \xi \vert^2 }+\vert B_0\cdot \xi \vert^2 }, \qquad \sigma_{2}(\xi)&=\Omega\frac {\xi_3}{\vert \xi \vert }-\sqrt{ \Omega^2\frac {\xi_3^2}{\vert \xi \vert^2 }+\vert B_0\cdot \xi \vert^2 },
\end{align}together with $\sigma_{3}(\xi)=-\sigma_{1}(\xi),\,\sigma_{4}(\xi)=-\sigma_{2}(\xi)$.
 \end{prop}
\begin{remark}
Here $\alpha_j,\beta_j:\mathbb R^3\to\mathbb C^3$, $j=1,\dots,4$,
satisfying
\begin{equation*}
\int_{\mathbb R^3}(1+|\xi|^2)^N
\bigl(|\alpha_j(\xi)|^2+|\beta_j(\xi)|^2\bigr)\,d\xi<\infty.   
\end{equation*}
\end{remark}

To study the dispersion, we restrict ourselves to the classical setting in \cite{lehnert1954magnetohydrodynamic} where the rotating axis and magnetic field are aligned $ B_0//e_3\,$:
\begin{cor}[linear decay]\label{linearcase}
Let $\alpha, \Omega\neq 0 $ and $B_0=\alpha e_3,$ and we define $\zeta = \frac {\alpha^2\vert \Omega\vert }{\Omega^2+\alpha^2} $. For $N\ge4$, let the initial data satisfy
\begin{equation*}
 u_{0},b_{0}\in
H^N(\R^3)
\cap\dot B^1_{1,1}(\R^3)
\cap\dot B^3_{1,1}(\R^3),   
\end{equation*}
the solution satisfies 
\begin{align*}
\Vert (u,b)(t)\Vert_{L^\infty}\lesssim \frac{1}{t\,\zeta}\Vert (u_{0},b_0)\Vert_{H^N
\cap\dot B_{1,1}^{1}\cap \dot B_{1,1}^{3} }.
\end{align*}
\end{cor}
We use the linear decay to extend the lifespan of a solution to the nonlinear equations. The resulting lifespan is quadratic up to a logarithmic correction.
Our main nonlinear result is the following: 
\begin{thm}\label{thm:extend}
Let $B_0=\alpha e_3$, with $\alpha,\Omega\neq 0$, and let $N>6$. Then there exist  constants $\varepsilon,c_0$ that any initial data $(u_{0},b_{0})$ satisfying
\begin{align*}
\|(u_{0},b_{0})\|_{H^N}=\varepsilon<c_0\min(1,\zeta)\, \text{ and }\,\|(u_{0},b_{0})\|_{W^{N-1,1}}\le \max \left(\eps, \frac{c_0\zeta}{\log(2+\zeta^2\varepsilon^{-2})}\right),
\end{align*}
we can obtain a unique solution $(u,b)\in C([0,T];H^N)$ to \eqref{eq:nl} such that
    \begin{align*}
\sup_{0\le t\le T}\|(u,b)\|_{H^N}(t)\le 2\varepsilon,
\end{align*}
and
\begin{align*}
T\ge
\frac{C\zeta}
{\varepsilon^2\log(2+\zeta^2\varepsilon^{-2})}
\end{align*}
for some $C>0$. 
\end{thm}
\subsection{Literature overview}
The nonlinear stability of Alfvén waves is established in  \cite{bardos1988longtime} for the ideal setting and in \cite{he2018global,wei2017global} for the dissipative setting. For a strictly monotone positive magnetic field $B=b(y)e_1$, the linearized dynamics exhibit damping \cite{ren2017linear}. The interplay of a monotone magnetic field and shear flow can further give rise to the formation of so-called magnetic islands in the linearized dynamics \cite{ren2021long,zhai2021long}.

A related line of research concerns the stability of the MHD equations near shear flows, especially near Couette flow, in both the ideal \cite{Knobelideal} and dissipative \cite{liss2020sobolev,Dolce,knobel2025sobolev,chen2024sobolev,zhao2024asymptotic,knobel2024nr,knobel2025suppression} settings. In that context, the dynamics are driven by shear-induced damping, current growth and enhanced dissipation, together with genuinely nonlinear effects such as echoes \cite{knobel2023echoes,knobel2025transition}.  For Taylor-Couette flow, the linear dynamics lead to exponential growth of the magnetic field due to the Ponamarenko dynamo \cite{navarro2025spectral}, yielding nonlinear instability \cite{navarro2025nonlinear}. These works rely on advection-induced mechanisms, whereas the present paper studies dispersive decay generated by the interplay between rotation and a background magnetic field. 

Dispersive effects have been studied for the Boussinesq
\cite{elgindi2015sharp,JurjaWidmayer2026},
rotating Euler
\cite{koh2014strichartz,ren2024global,guo2023global},
rotating Navier--Stokes \cite{ko2024global},
beta-plane \cite{pusateri2018global},
and Euler--Maxwell \cite{guo2016global} equations.
For rapidly rotating MHD with vanishing viscosity,
Charve and Yao~\cite{CharveYao} used Strichartz and
energy estimates to obtain global strong solutions
and convergence rates.
The rotating Euler and stratified Boussinesq models
feature homogeneous phases such as $\frac{\xi_3}{|\xi|}$
and $\frac{\xi_1}{|\xi|}$.
Nonhomogeneous radial phases were studied by
Melinand~\cite{Melinand2024}, with applications
to water wave models.
Here, we consider phases $\xi_3h(|\xi|)$,
which need not be homogeneous, and apply their
dispersive estimates to derive lifespan bounds
for rotating ideal MHD.In contrast to these homogeneous dispersion relations, magneto-Coriolis waves arise from the coupling of two distinct restoring mechanisms, namely rotation and magnetic tension, and possess several dispersion branches with nonhomogeneous radial dependence. The main feature of our approach is to formulate the dispersive analysis through verifiable radial and angular nondegeneracy conditions for phases of the form $\sigma(\xi)=\xi_3 h(|\xi|).$
This provides a unified treatment of the different magneto-Coriolis branches and their frequency regimes.

Magneto-Coriolis waves arise on length scales where the background magnetic field may be treated as constant. On larger scales, the variation of the global magnetic field becomes relevant and at leading order yields a balance between the Lorentz force and rotation, described by the Taylor equations \cite{taylor1963magneto,ferriz2019advances}. Even in that regime, well-posedness is a major challenge and is known only in special settings \cite{gallagher2017wellposedness}. Heuristically, low frequencies are more strongly influenced by the large-scale magnetic field, whereas higher frequencies and the fast magneto-Coriolis branch retain dispersive behavior.

\subsection{Radius dependent dispersive estimates}
To simplify some expressions, we denote 
\begin{align*}
r&=r(\xi)=\vert \xi\vert,\\
\theta &=\theta(\xi)=\frac {\xi_3}r ,\\
\gamma&=\tfrac 13 \left(1-\tfrac {r h''(r)}{h'(r)}\right).
\end{align*}

To obtain the dispersion of magneto-Coriolis waves, we consider the more general dispersion relations of the form 
\begin{align}
    \sigma (\xi )= \xi_3 h(r),
\end{align}
and establish dispersive estimates for 
\begin{align}\label{dispersive}
    I_h[f](t,x)&= \int \exp(i (\xi \cdot x - \sigma(\xi) t )) \hat f(\xi) d\xi.
\end{align}
For general $h$, $I_h$ will not decay in the supremum norm. A simple example is when $h=1$. Therefore, we require sufficient conditions on $h$ such that $I_h$ decays in time:

\begin{ass}\label{ass:h}
Let $h\in C^4((0,\infty);\R)$.
Assume that \textnormal{(A)}--\textnormal{(D)} below hold
for some constants $c,a_0,A_0>0$.
Here $\chi\in C^\infty([-1,1])$ is fixed,
with 
\begin{equation*}
 0\notin\supp\chi,
\qquad
A_0 s^2\le\frac14
\quad\text{for every } s \in \supp(1-\chi).   
\end{equation*}

\begin{enumerate}
\item[\textnormal{(A)}]
There exist $m_0,m_\infty\ge0$ and $\zeta_0,\zeta_\infty>0$
such that, for every $r>0$,
\begin{equation*}
  |h'(r)|
\ge c\min\left\{\zeta_0r^{-m_0},
                 \zeta_\infty r^{-m_\infty}\right\};
\end{equation*}

\item[\textnormal{(B)}]
The number of stationary points satisfies
\begin{equation*}
 \sup_{x\in\R^3,\ n\in\mathbb Z}
\#\left\{
\xi:
2^n\le|\xi|\le2^{n+1},\
\chi\left(\frac{\xi_3}{|\xi|}\right)\ne0,\
\nabla_\xi\bigl(x\cdot\xi-\xi_3h(|\xi|)\bigr)=0
\right\}
\le c^{-1};  
\end{equation*}

\item[\textnormal{(C)}]
For every $r>0$,
\begin{equation*}
0<a_0
\le 4-3\gamma(r)=3+\tfrac {r h''(r)}{h'(r)}
\le A_0;
\end{equation*}

\item[\textnormal{(D)}]
For every $r>0$,
\begin{equation*}
  r^2|h'''(r)|
\le c^{-1}|h'(r)|.  
\end{equation*}
\end{enumerate}
\end{ass}

\begin{remark}
The proof of Theorem~\ref{prop:strich} below uses
conditions (A), (C), and (D).
We write (B) as an additional geometric property and
verify it for the magneto--Coriolis profiles in
Lemma~\ref{lem:ms}.
\end{remark}

\begin{thm}\label{prop:strich}
Let $f\in C^\infty_c(\R^3)$ and let $h$ satisfy
Assumption~\ref{ass:h}. Then, for any  $t>0$,
\begin{equation*}
 \|I_h[f](t)\|_{L^\infty}
\le
\frac{C}{t}
\left(
\zeta_0^{-1}\|f\|_{\dot B_{1,1}^{1+m_0}}
+
\zeta_\infty^{-1}\|f\|_{\dot B_{1,1}^{1+m_\infty}}
\right).   
\end{equation*}
In particular, if $\zeta_0=\zeta_\infty=\zeta$, then we have
\begin{equation*}
\|I_h[f](t)\|_{L^\infty}
\le
\frac{C}{t\,\zeta}
\left(
\|f\|_{\dot B_{1,1}^{1+m_0}\cap
\dot B_{1,1}^{1+m_\infty}}
\right).    
\end{equation*}
\end{thm}
\begin{remark}
The $t^{-1}$ decay rate is optimal.
Indeed, since $h'(r)\neq0$, the Hessian has rank two
on the equatorial plane $\xi_3=0$. The $|t|^{-3/2}$ decay associated with a fully
nondegenerate three-dimensional Hessian cannot be
obtained by applying stationary phase in two
nondegenerate directions.
For example, $h(r)=r$ gives the phase $\xi_3|\xi|$
of the three-dimensional Benjamin--Ono equation,
up to a change of sign and a permutation of coordinates.
For this phase, the frequency-localized $|t|^{-1}$
decay estimate is optimal; see
\cite[Proposition~3.1]{hickman2019higher}.
The choice $h(r)=r^{-1}$ corresponds to inertial waves;
see \cite{koh2014strichartz} for the associated
Strichartz estimates.
\end{remark}

\subsection{Proof strategy and overview}
Inspired by the dispersive approach to lifespan
estimates in \cite{elgindi2015sharp},
we combine linear decay with Sobolev energy estimates.
The argument has two main components. In Section~\ref{sec:disperse}, we prove an abstract dispersive
estimate for phases of the form $\xi_3 h(|\xi|)$.
Using an explicit computation of the Hessian and a fixed
smooth angular decomposition, we identify two nondegenerate
directions on each localized piece.
Then we apply two-dimensional stationary phase to obtain
the $t^{-1}$ decay estimate. In Section \ref{sec:mcw}, we derive the explicit representation of solutions to the linearized equations and verify that the corresponding dispersion relations fit the abstract framework. This yields the linear decay estimate, which is then combined with Sobolev energy estimates to obtain the extended lifespan of the nonlinear equations.

\section{Preliminaries}
\begin{definition}
 For an integer $N\ge0$ and $1\le p\le\infty$, we define
\begin{equation*}
W^{N,p}(\mathbb R^3)
=
\left\{
f\in L^p(\mathbb R^3):
\partial^\beta f\in L^p(\mathbb R^3)
\text{for all }|\beta|\le N
\right\},    
\end{equation*}
with norm
\begin{equation*}
\begin{cases}
\|f\|_{W^{N,p}}
=\sum\limits_{|\beta|\le N}\|\partial^\beta f\|_{L^p},\\   \|f\|_{W^{N,\infty}}
=
\sum\limits_{|\beta|\le N}\|\partial^\beta f\|_{L^\infty}
    \end{cases}
\end{equation*}

\end{definition}
For $N\ge0$, we have $W^{N,2}=H^N$,
with equivalent norms.
\begin{definition}\label{inhomosob}
Let $N \in \R$. The Sobolev space $H^N(\R^d)$ consists of tempered distributions $u$ whose Fourier transform belongs to $L^2_{\mathrm{loc}}(\R^d)$ and satisfies
\begin{equation*}
 \|u\|_{H^N}^2 
\stackrel{\mathrm{def}}{=} 
\int_{\R^d} (1 + |\xi|^2)^N |\hat{u}(\xi)|^2 \, d\xi < \infty.  \end{equation*}
\end{definition}
\begin{prop}\label{HLP}
Let $\varphi,\psi\in C_c^\infty(\R^3\setminus\{0\})$ be radial
functions such that
\begin{equation*}
\supp\varphi
\subset
\left\{\xi\in\R^3:\frac12\leq|\xi|\leq2\right\},
\qquad
\sum_{j\in\Z}\varphi(2^{-j}\xi)=1,
\quad \xi\neq0,    
\end{equation*}

and
\begin{equation*}
\psi=1\quad\text{on }\supp\varphi,
\qquad
\supp\psi
\subset
\left\{\xi\in\R^3:\frac14\leq|\xi|\leq4\right\}.    
\end{equation*}

For $j\in\Z$, define
\begin{equation*}
\dot\Delta_j f
:=
\mathcal F^{-1}
\left(
\varphi(2^{-j}\xi)\widehat f(\xi)
\right)    
\end{equation*}
and
\begin{equation*}
 S_j f
:=
\mathcal F^{-1}
\left(
\psi(2^{-j}\xi)\widehat f(\xi)
\right).   
\end{equation*}
Then
\begin{equation*}
\dot\Delta_j f=\dot\Delta_j(S_j f).    
\end{equation*}
Moreover, for every $s\in\R$ and $1\leq p\leq\infty$,
\begin{equation*}
 \sum_{j\in\Z}
2^{js}\|S_j f\|_{L^p}
\lesssim
\sum_{j\in\Z}
2^{js}\|\dot\Delta_j f\|_{L^p}.   
\end{equation*}
\end{prop}

\begin{definition}
Let $s\in\R$ and $1\leq p,q\leq\infty$. The homogeneous Besov
space $\dot B^s_{p,q}(\R^3)$ consists of all
$f\in\mathcal S'(\R^3)$, such that
\begin{equation*}
\|u\|_{\dot B^s_{p,q}}
:=
\left(
\sum_{j\in\Z}
2^{jsq}\|\dot\Delta_j u\|_{L^p}^q
\right)^{1/q}
<\infty.    
\end{equation*}
When $q=\infty$, we use the usual modification
\begin{equation*}
\|f\|_{\dot B^s_{p,\infty}}
:=
\sup_{j\in\Z}
2^{js}\|\dot\Delta_j f\|_{L^p}.  
\end{equation*}
\end{definition}
\begin{prop}\label{embedding}
For every positive integer $N$, $1\le p\le\infty$,
and $0<s<N$, we have the continuous embedding
\begin{equation*}
W^{N,p}(\mathbb R^3)
\hookrightarrow
\dot B^s_{p,1}(\mathbb R^3),    
\end{equation*}
with
\begin{equation*}
\|f\|_{\dot B^s_{p,1}}
\le C\|f\|_{W^{N,p}}.   \end{equation*}

\end{prop}
\begin{prop}
\label{sobolev-product}
Let $N>\frac32$. There exists a constant $C_N>0$
such that, for all $f,g\in H^N(\mathbb R^3)$,
the product $fg$ belongs to $H^N(\mathbb R^3)$ and
\begin{equation*}
\|fg\|_{H^N}
\le C_N\|f\|_{H^N}\|g\|_{H^N}.    
\end{equation*}
\end{prop}
\section{Proof of dispersive estimates}\label{sec:disperse}

We prove Theorem \ref{prop:strich} by the stationary phase approximation and the explicit representation of the eigenvalues. 

The first lemma classifies the eigenvalues of the Hessian matrix:
\begin{lemma}\label{lem:hess}
For $h$ such that $h'(r)\neq 0$ and let  $\sigma\in C^2(\R^3\setminus \{0\}), \  \sigma(\xi ) = \xi_3 h(r)$ then the Hessian matrix $H\sigma(\xi) $ has the eigenvalues 
    \begin{align*}
    \lambda_0(\xi) &= \theta h'(r),\\
    \lambda_\pm(\xi)  &= \frac 1 2h' \left(\theta (4-3\gamma )  \mp \sqrt{(3 \theta \gamma)^2+4(1-3 \theta^2 \gamma)}\right),
    \end{align*}
    for almost all $\xi\in \R^3$. 
\end{lemma}
\textbf{{Proof of Theorem \ref{prop:strich}}:}
By Lemma \ref{lem:hess}, we have eigenvalues $\lambda_0, \lambda_\pm $ that satisfy 
\begin{equation}\label{lamda12}
\begin{aligned}
 \lambda_+(\xi)\lambda_-(\xi)&= \frac 1 4(h')^2 \left(\theta^2 (4-3\gamma )^2  - (3 \theta \gamma)^2-4(1-3 \theta^2 \gamma)\right)\\
    &= -(h')^2\left(1-\theta^2 (4-3\gamma)\right),    
\end{aligned}
\end{equation}
and so 
\begin{align}\label{Hmatrix}
\det H\sigma &=\lambda_0\lambda_+\lambda_- = -(h')^3\left(1-\theta^2 (4-3\gamma)\right)\theta. 
\end{align}
Thus, the Hessian matrix is degenerate when $\theta=0$.
And we will identify two uniformly nondegenerate
directions on each part of the angular decomposition.

Let $\varphi$ and $S_n$ be as in Proposition~\ref{HLP},
and set $\varphi_n(\xi)=\varphi(2^{-n}\xi)$. Using $\dot\Delta_n f=\dot\Delta_n(S_n f)$
from Proposition \ref{HLP} and Young inequality, we obtain
\begin{align}\label{In}
\|I_h[f](t)\|_{L^\infty}
&=\left\|I_h\left[\sum_{n\in\mathbb Z}\dot\Delta_n f\right](t)
\right\|_{L^\infty}\\
&\lesssim
\sum_{n\in\mathbb Z}
\left\|2^{(k-3)n}I_h[\mathcal F^{-1}\varphi_n]\right\|_{L^\infty}
\,2^{(3-k)n}\|S_n f\|_{L^1}
\end{align}
where $k$ will be chosen below.

Then we establish a estimate on
\begin{align}
   T_n(x):=2^{(k-3)n} \int  \exp(i (\xi \cdot x - \sigma(\xi) t )) \varphi(2^{-n} \xi) d\xi.
\end{align}
By changing the variables $\eta =2^{-n} \xi$ and $\tilde x= 2^n x$, we obtain  
\begin{align*}
    T_n (\tilde x) &= 2^{kn}\int d\eta  \exp(i (\eta \cdot \tilde x - \tilde \sigma_{n}(\eta ) t )) \varphi(\eta ), \\
    &= 2^{kn}\int d\eta  \exp(i (\eta \cdot \tilde x - \tilde \sigma_{n }(\eta ) t )) \varphi(\eta ),
\end{align*}
with $\tilde \sigma_{n}(\eta )=  \sigma(2^n\eta )$.
We decompose
\begin{equation*}
 T_n=T_{n,1}+T_{n,2},   
\end{equation*}
and use the cutoff $\chi$ fixed in
Assumption~\ref{ass:h}.
Since $0\notin\supp\chi$, we may choose $\delta>0$
such that
\begin{equation*}
\supp\chi\subset\{s\in[-1,1]:|s|\ge\delta\}. 
\end{equation*}
So we write
\begin{align*}
T_{n,1}(\tilde x)
&=
2^{kn}\int_{\R^3}
e^{i(\eta\cdot\tilde x-t\tilde\sigma_n(\eta))}
\varphi(\eta)
\left(1-\chi\left(\frac{\eta_3}{|\eta|}\right)\right)
\,d\eta,\\
T_{n,2}(\tilde x)
&=
2^{kn}\int_{\R^3}
e^{i(\eta\cdot\tilde x-t\tilde\sigma_n(\eta))}
\varphi(\eta)
\chi\left(\frac{\eta_3}{|\eta|}\right)
\,d\eta.
\end{align*}

\textbf{Bound on $T_{n,1}$.}

On the support of
\begin{equation}\label{supp}
\varphi(\eta)
\left(1-\chi\left(\frac{\eta_3}{|\eta|}\right)\right), \end{equation}we have $\theta=\frac{\eta_3}{|\eta|}\in\supp(1-\chi)$.
Then by the 
Assumption~\ref{ass:h}(C), we know that
\begin{equation*}
\theta^2(4-3\gamma(r))
\le A_0\theta^2\le\frac14,
\qquad \text{where}\qquad r=2^n|\eta|.    
\end{equation*}
Then from \eqref{lamda12} we know that
\begin{equation}\label{eq:two-eigenvalues}
|\lambda_+(2^n\eta)\lambda_-(2^n\eta)|
\ge \frac34|h'(2^n|\eta|)|^2.
\end{equation}
Since
$\xi=2^n\eta$ satisfies
\begin{equation*}
2^{n-1}\le|\xi|=2^n|\eta|\le2^{n+1}.    
\end{equation*}
We compare $|h'(2^n|\eta|)|$
with $|h'(2^n)|$,  we use Assumptions~\ref{ass:h}(A) and (C):
\begin{equation*}
\left|\frac{rh''(r)}{h'(r)}\right|
\le M,
\qquad
M=\max\{|a_0-3|,|A_0-3|\}.    
\end{equation*}
Thus, for $\frac14\le|\eta|\le4$,
\begin{align*}
\left|
\log\frac{|h'(2^n|\eta|)|}{|h'(2^n)|}
\right|
&=
\left|
\int_{2^n}^{2^n|\eta|}
\frac{h''(s)}{h'(s)}\,ds
\right|\\
&\le M|\log|\eta||
\le M\log4.
\end{align*}
So we have
\begin{equation}\label{eq:hprime-comparison}
m_0
\le
\frac{|h'(2^n|\eta|)|}{|h'(2^n)|}
\le m.
\end{equation}

By the chain rule,
\begin{equation}\label{eq:chain}
H_\eta\widetilde\sigma_n(\eta)
=
2^{2n}H_\xi\sigma(\xi).   
\end{equation}
Let $p_+,p_-$ be orthonormal eigenvectors of
$H_\xi\sigma(\xi)$ corresponding to $\lambda_+,\lambda_-$. Then
\begin{equation*}
   B_n(\eta)= \begin{pmatrix}
p_+^T H_\eta\widetilde\sigma_n(\eta)p_+
&
p_+^T H_\eta\widetilde\sigma_n(\eta)p_-\\
p_-^T H_\eta\widetilde\sigma_n(\eta)p_+
&
p_-^T H_\eta\widetilde\sigma_n(\eta)p_-
\end{pmatrix}
=
2^{2n}
\begin{pmatrix}
\lambda_+(\xi)&0\\
0&\lambda_-(\xi)
\end{pmatrix}.
\end{equation*}
Combining \eqref{eq:two-eigenvalues} and \eqref{eq:hprime-comparison}, we know that in \eqref{supp} we have:
\begin{equation*}
  \left|
\det\left(
2^{2n}
\begin{pmatrix}
\lambda_+(\xi)&0\\
0&\lambda_-(\xi)
\end{pmatrix}
\right)
\right|
=
2^{4n}|\lambda_+(\xi)\lambda_-(\xi)|
\ge
\frac34\,2^{4n}|h'(2^n|\eta|)|^2
\ge
m_0\,2^{4n}|h'(2^n)|^2.  
\end{equation*}

First, we know that Assumptions~\ref{ass:h}(C) and (D) give
\begin{equation*}
   |h''(r)|\lesssim \frac{|h'(r)|}{r},
\qquad
|h'''(r)|\lesssim  \frac{|h'(r)|}{r^2}. 
\end{equation*}Then we differentiate $\sigma(\xi)=\xi_3h(|\xi|)$, we can obtain
\begin{equation*}
  \|D_\xi^2\sigma(\xi)\|
\le C\bigl(|h'(r)|+r|h''(r)|\bigr)
\le C|h'(r)|,  
\end{equation*}
and
\begin{equation*}
\|D_\xi^3\sigma(\xi)\|
\lesssim \left(
\frac{|h'(r)|}{r}
+|h''(r)|
+r|h'''(r)|
\right)
\lesssim \frac{|h'(r)|}{r}.
\end{equation*}
Using \eqref{eq:chain} and \eqref{eq:hprime-comparison},
we obtain, for $\frac14\le|\eta|\le4$,
\begin{equation}\label{eq:phase-derivative-bounds}
\|D_\eta^2\widetilde\sigma_n(\eta)\|
+
\|D_\eta^3\widetilde\sigma_n(\eta)\|
\lesssim \,2^{2n}|h'(2^n)|.
\end{equation}

Fix a point $\eta^0$ in \eqref{supp}.
Choose $p_+,p_-$ at $\xi^0=2^n\eta^0$ as above,
and keep these two vectors fixed. We write
\begin{equation*}
 |\det B_n(\eta^0)|
\ge m_0\,2^{4n}|h'(2^n)|^2.   
\end{equation*}
Then we take $|\eta-\eta^0|<\rho$, without loss of generality where we can choose $\rho\le\frac18$,
so that the segment joining $\eta^0$ to $\eta$
lies in $\{\frac14\le|\eta|\le4\}$.
Since $p_+,p_-$ are fixed, the mean value theorem
and \eqref{eq:phase-derivative-bounds} give
\begin{equation*}
|\det B_n(\eta)-\det B_n(\eta^0)|
\le
C_1\,2^{4n}|h'(2^n)|^2\,|\eta-\eta^0|.    
\end{equation*}
Choose
\begin{equation*}
 0<\rho\le
\min\left\{\frac18,\frac{m_0}{2C_1}\right\}.   
\end{equation*}
Then, for every $\eta\in B(\eta^0,\rho)$,
\begin{align*}
|\det B_n(\eta)|
&\ge
|\det B_n(\eta^0)|
-
|\det B_n(\eta)-\det B_n(\eta^0)|\\
&\ge
(m_0-C_1\rho)\,2^{4n}|h'(2^n)|^2\\
&\ge
\frac{m_0}{2}\,2^{4n}|h'(2^n)|^2.
\end{align*}
The radius $\rho$ is independent of $n$ and of
the center $\eta^0$.
Since the support of the amplitude is a fixed compact set,
we can cover it by finitely many balls
$B(\eta^\ell,\rho/4)$, $\ell=1,\ldots,N$.
Choose functions $a_\ell$ such that
\begin{equation*}
 \varphi(\eta)
\left(1-\chi\left(\frac{\eta_3}{|\eta|}\right)\right)
=
\sum_{\ell=1}^N a_\ell(\eta),
\qquad
\supp a_\ell\subset B(\eta^\ell,\rho/2).   
\end{equation*}

At each center $\eta^\ell$, choose the two fixed
directions $p_\ell,q_\ell$ as above, and complete them
to an orthonormal basis $(p_\ell,q_\ell,w_\ell)$.
And we write
\begin{equation*}
\eta_\ell(u,v)
=
\eta^\ell+u_1p_\ell+u_2q_\ell+vw_\ell,
\qquad u=(u_1,u_2).    
\end{equation*}
Since this change of coordinates is orthogonal.
Moreover, the support condition implies
$|u|^2+v^2<\rho^2/4$.
Thus,
\begin{align*}
T_{n,1}(\widetilde x)
=
2^{kn}\sum_{\ell=1}^N
\int_{-\rho/2}^{\rho/2}
\int_{\mathbb R^2}
&e^{i\left(
\eta_\ell(u,v)\cdot\widetilde x
-t\widetilde\sigma_n(\eta_\ell(u,v))
\right)}\\
&\qquad\qquad
\times a_\ell(\eta_\ell(u,v))\,du\,dv.
\end{align*}

For each fixed $v$, the Hessian in $u$ satisfies
\begin{equation*}
D_u^2\left[
\widetilde\sigma_n(\eta_\ell(u,v))
\right]
=
B_{n,\ell}(\eta_\ell(u,v)),  
\end{equation*}
where $B_{n,\ell}$ is the matrix
similar as $B_n$.
Hence,
\begin{equation*}
\left|
\det D_u^2\left[
\widetilde\sigma_n(\eta_\ell(u,v))
\right]
\right|
\ge m_0\,2^{4n}|h'(2^n)|^2.    
\end{equation*}
The second and third derivatives in $u$ are bounded by
$C2^{2n}|h'(2^n)|$.

We observe that $\left(
\eta_\ell(u,v)\cdot\widetilde x
-t\widetilde\sigma_n(\eta_\ell(u,v))
\right)$, then $t2^{2n}|h'(2^n)|\ge1$,
the two-dimensional stationary-phase estimate gives
\begin{align*}
\left|
\int_{\mathbb R^2}
e^{i\left(
\eta_\ell(u,v)\cdot\widetilde x
-t\widetilde\sigma_n(\eta_\ell(u,v))
\right)}
a_\ell(\eta_\ell(u,v))\,du
\right|
\le
\frac{C}{t2^{2n}|h'(2^n)|}.
\end{align*}

Integrating in $v$ and summing over $\ell$, we obtain
\begin{equation*}
 |T_{n,1}(\widetilde x)|
\le
2^{kn}\sum_{\ell=1}^N
\int_{-\rho/2}^{\rho/2}
\frac{C}{t2^{2n}|h'(2^n)|}\,dv
\le
\frac{C2^{kn}}{t2^{2n}|h'(2^n)|}.   
\end{equation*}
On the other hand, we know that
\begin{equation*}
|T_{n,1}(\widetilde x)|
\le
2^{kn}
\int_{\mathbb R^3}
\left|
\varphi(\eta)
\left(1-\chi\left(\frac{\eta_3}{|\eta|}\right)\right)
\right|\,d\eta
\le C2^{kn}.    
\end{equation*}
Combining these two bounds yields
\begin{equation}\label{eq:Tn1}
\|T_{n,1}\|_{L^\infty}
\le
C2^{kn}
\min\left\{
1,\frac{1}{t2^{2n}|h'(2^n)|}
\right\},
\qquad t>0.
\end{equation}

\textbf{Bound on $T_{n,2}$.}

Since $0\notin\supp\chi$, there exists $\delta>0$ such that
\begin{equation*}
|\theta|\ge\delta
\quad\text{on }\supp\chi.    
\end{equation*}
For vectors $p,q\in\xi^\perp$, the Hessian satisfies
\begin{equation*}
p^T H_\xi\sigma(\xi)q=\theta h'(r)\,p\cdot q.    
\end{equation*}
Thus, in any orthonormal basis of $\xi^\perp$,
the restriction of $H_\xi\sigma$ is
$\theta h'(r)\mathrm{Id}_2$.
Then we obtain
\begin{equation*}
 \theta^2
\frac{|h'(2^n|\eta|)|^2}{|h'(2^n)|^2}
\ge
\delta^2
\frac{|h'(2^n|\eta|)|^2}{|h'(2^n)|^2}
\gtrsim1.   
\end{equation*}

Repeating the preceding localization, applying
two-dimensional stationary phase in $(u_1,u_2)$,
and integrating in $v$, we obtain
\begin{equation}\label{eq:Tn2}
\|T_{n,2}\|_{L^\infty}
\lesssim
2^{kn}
\min\left\{
1,\frac{1}{t2^{2n}|h'(2^n)|}
\right\}.
\end{equation}

\medskip
Combining \eqref{eq:Tn1} and \eqref{eq:Tn2},
and using Assumption~\ref{ass:h}(A), 
\begin{equation*}
\frac{2^n}{|h'(2^n)|}
\lesssim
\zeta_0^{-1}2^{(1+m_0)n}+
\zeta_\infty^{-1}2^{(1+m_\infty)n}, n\ge1.   
\end{equation*}
We obtain
\begin{align*}
\|I_h[\mathcal F^{-1}\varphi_n](t)\|_{L^\infty}
&=2^{(3-k)n}\|T_n\|_{L^\infty}\\
&\lesssim
2^{3n}
\min\left\{
1,\frac{1}{t2^{2n}|h'(2^n)|}
\right\}\\
&\le\frac{C}{t}
\frac{2^n}{|h'(2^n)|}\\
&\lesssim
\frac{1}{t}
\left(
\zeta_0^{-1}2^{(1+m_0)n}
+
\zeta_\infty^{-1}2^{(1+m_\infty)n}
\right).
\end{align*}
Hence, Young inequality and
Proposition~\ref{HLP} give
\begin{align*}
\|I_h[f](t)\|_{L^\infty}
&\lesssim
\frac{1}{t}
\sum_n
\left(
\zeta_0^{-1}2^{(1+m_0)n}
+
\zeta_\infty^{-1}2^{(1+m_\infty)n}
\right)
\|S_nf\|_{L^1}\\
&\lesssim
\frac{1}{t}
\left(
\zeta_0^{-1}\|f\|_{\dot B_{1,1}^{1+m_0}}
+
\zeta_\infty^{-1}
\|f\|_{\dot B_{1,1}^{1+m_\infty}}
\right).
\end{align*}

This proves the desired estimate.

\subsection{Eigenvalues of the Hessian matrix}
In this section, we prove Lemma \ref{lem:hess}, without loss of generality, we assume that $\xi_3\ge 0$ and $h'(r)<0$. We calculate 
\begin{align}\label{de1}
    \p_i \sigma &= \delta_{i3} h(r) +\frac{\xi_3 \xi_i}r h'(r)
\end{align}
and 
\begin{align*}
    \p_j \p_i \sigma &=\left( \delta_{i3}\frac{\xi_j}{r} +\delta_{j3}\frac{\xi_i}{r}  \right)h'(r) +\frac{\xi_3}r \left(\delta_{ij}h'(r)+\frac{\xi_i\xi_j }{r^2}(r h''(r)-h'(r) \right)\\
    &= \left(\delta_{ij}-\frac{\xi_i\xi_j }{r^2}\right)\frac{\xi_3}rh'(r) +\frac{\xi_i\xi_j }{r^2}\xi_3 h''(r)+\left( \delta_{i3}\frac{\xi_j}{r} +\delta_{j3}\frac{\xi_i}{r}  \right)h'(r).
\end{align*}
and so we denote
\begin{align*}
    H\sigma &=  H_1+H_2+H_3,\\
    H_1&= (I-\tfrac{\xi\xi^T}{r^2})\cdot \frac{\xi_3}rh'(r),\\
    H_2&=\tfrac{\xi\xi^T}{r^2}\cdot\xi_3 h''(r),\\
    H_3&=\left( \delta_{i3}\frac{\xi_j}{r} +\delta_{j3}\frac{\xi_i}{r}  \right)h'(r)\\
    &= \begin{pmatrix}
        0&0&\xi_1\\
        0&0&\xi_2\\
        \xi_1&\xi_2&2\xi_3
    \end{pmatrix}\frac{h'(r)}r.
\end{align*}
For 
\begin{align*}
    v_0&= \begin{pmatrix}
        -\xi_2\\
        \xi_1\\ 
        0
    \end{pmatrix}, 
\end{align*}
we obtain $H_2v_0=H_3v_0=0$ and 
\begin{align*}
Hv_0=H_1v_0 = \frac{\xi_3}rh'(r)v_0,
\end{align*}
and so we have the eigenvalue $\lambda_0 = \frac{\xi_3}rh'(r)$.

To calculate the eigenvectors $w_\pm$, we observe the structure of $H_3$ and consider the eigenvectors of $H_3$, 
\begin{align}
    w_\pm &= \xi \pm e_3 \vert \xi\vert 
\end{align}
which satisfy 
\begin{equation}\label{H3}
H_3w_\pm=\frac{h'}r\left(\xi_3\pm \vert \xi \vert\right)w_\pm.    
\end{equation}
To obtain the eigenvectors of $H$ define 
\begin{equation*}
    \tilde w=\tilde w[g] =  \xi + e_3 \vert \xi\vert g(\xi),
\end{equation*}
and we determine two functions $g(\xi)\neq 0 $ such that $\tilde w$ is a eigenvector. We obtain 
\begin{equation*}
    H_1\tilde w =a(\xi)\left(\tilde w -\xi \left(1+g(\xi)\frac {\xi_3}r\right)\right),  \end{equation*}
where $a(\xi)=\frac{\xi_3}rh'(r),$ and
\begin{equation*}
    H_2\tilde w=b(\xi)\xi \left(1+g(\xi)\frac {\xi_3}r\right),
\end{equation*}
    where $b(\xi)=\xi_3 h''(r).$
Since $\tilde w= \frac{1+g}2w_++\frac{1-g}2w_-$, we can obtain 
\begin{equation*}
H_3\tilde w =\frac{h'}r (\xi_3+\tfrac r g)\tilde w +h'(g-\tfrac 1 g )\xi.
\end{equation*}
Therefore
\begin{equation*}
    H\tilde w = (a+\xi_3\tfrac{h'}r+\tfrac {h'} g  )\tilde w +\left((b-a) (1+g\tfrac {\xi_3}r)+ h'(g-\tfrac 1 g )  \right)\xi.
\end{equation*}
In particular,  we have an eigenvalue for $g$ such that 
\begin{equation*}
    \left((b-a) (1+g\tfrac {\xi_3}r)+ h'(g-\tfrac 1 g )  \right)=0,
\end{equation*}
and so 
\begin{equation*}
    g^2 (h'+(b-a)\tfrac {\xi_3}r)+(b-a) g-h' =0.
\end{equation*}
This is solved by
\begin{equation*}
 g_\pm(\xi) = \frac{a-b \pm \sqrt{(a-b)^2+4h'(h'+(b-a)\tfrac {\xi_3}r) }}{2\left(h'+(b-a)\tfrac {\xi_3}r\right)}. 
\end{equation*}

Since $h'\neq 0$, the $g_\pm$ is well defined and non zero if $0 \neq h'+(b-a)\tfrac {\xi_3}r$. Since the set $\xi$ such that $0 = h'+(b-a)\tfrac {\xi_3}r$ has measure $0$ it can be neglected and it is sufficient to consider $\xi$ such that $g_\pm(\xi)$ is well defined and non zero. Then  $v_\pm =\tilde w[g_\pm ]$ are eigenvectors and satisfy
\begin{equation*}
    H v_\pm =\lambda_\pm v_{\pm} 
\end{equation*}
with the corresponding eigenvalue
\begin{align*}
    \lambda_\pm &= a+\xi_3\frac{h'}r+\frac{h'} {g_\pm  } \\
    &= 2\frac{\xi_3}r h'+ 2h'\frac{h'+(b-a)\tfrac {\xi_3}r}{a-b \pm \sqrt{(a-b)^2+4h'(h'+(b-a)\tfrac {\xi_3}r) }}.
\end{align*}
Now we simplify the eigenvalues. Note that  $g$ and $\lambda$ are only dependent on the variables $r$ and $\theta=\frac {\xi_3} r$. We have
\begin{equation*}
  g_\pm(r,\theta)= \frac{a-b \pm \sqrt{(a-b)^2+4h'(h'+(b-a)\theta) }}{2\left(h'+(b-a)\theta \right)},   
\end{equation*}
and by inserting $a(\xi)=h'\theta$ we obtain 
\begin{align*}
    \lambda_\pm(r,\theta) &= 2h' \left(\theta + \frac{h'+(b-a)\theta}{a-b \pm \sqrt{(a-b)^2+4h'(h'+(b-a)\theta) }}\right).
\end{align*}
We use that $a-b= \theta (h'-r h'')=3 \theta h'\underbrace{(\frac 1 3 -\frac{r h''}{3h'})}_{=\gamma(r)} $
\begin{align*}
    \lambda_\pm(r,\theta) &= 2h' \left(\theta + \frac{h'-3 \theta h'\gamma\theta}{3 \theta h'\gamma \pm \sqrt{(3 \theta h'\gamma)^2+4h'(h'-3  h'\theta^2\gamma) }}\right).
\end{align*}
Since  $h'(r)< 0$
\begin{equation*}
    \begin{aligned}
     \lambda_\pm (r,\theta)&= 2h' \left(\theta + \frac{1-3 \theta^2 \gamma}{3 \theta \gamma \mp \sqrt{(3 \theta \gamma)^2+4(1-3 \theta^2 \gamma) }}\right)\\
    &= 2h'\left(\theta -\frac{1-3 \theta^2 \gamma}{4(1-3 \theta^2 \gamma) }\left(3 \theta \gamma \pm \sqrt{(3 \theta \gamma)^2+4(1-3 \theta^2 \gamma)}\right)\right)\\
    &= 2h' \left(\theta -\frac{1}{4 }\left(3 \theta \gamma \pm \sqrt{(3 \theta \gamma)^2+4(1-3 \theta^2 \gamma)}\right)\right)\\
    &= 2h' \left(\theta (1-\frac{3}{4 }\gamma )  \mp \frac 14 \sqrt{(3 \theta \gamma)^2+4(1-3 \theta^2 \gamma)}\right)\\
    &= \frac 1 2h' \left(\theta (4-3\gamma )  \mp \sqrt{(3 \theta \gamma)^2+4(1-3 \theta^2 \gamma)}\right).   
    \end{aligned}
\end{equation*}
Thus, we obtain Lemma \ref{lem:hess}.

The second lemma is a technical estimate used to bound the degenerate region of the eigenvalues: 
\begin{lemma}\label{lem:assd}
Let Assumption \ref{ass:h} C holds, then we have 
        $$\sup_n \vert \{\eta: \ \left\vert 1- \tfrac {\eta_3^2}{\vert  \eta \vert^2}(4-3\gamma(2^n\vert \eta\vert))\right\vert\le \lambda^2 \}\cap \{2^{-1}\le \vert \eta \vert \le 2\}\vert\lesssim \lambda^2 .$$
\end{lemma}
\begin{proof}
Without loss of generality, we assume $0<\lambda\le 1/10$, since for
$\lambda>1/10$ the desired estimate follows from the finite
measure of the annulus $\{2^{-1}\le|\eta|\le2\}$.
We fix $n$, write $\gamma_n =\gamma(2^n\vert \eta\vert ) $ and denote 
\begin{align*}
\Gamma_n &=\{\eta: \ \left\vert 1- \tfrac {\eta_3^2}{\vert  \eta \vert^2}(4-3\gamma_n)\right\vert\le \lambda^2 \}\\&=\{\eta: \ 1-\lambda^2\le  \tfrac {\eta_3^2}{\vert  \eta \vert^2}(4-3\gamma_n)\le 1+\lambda^2 \}.
\end{align*}
We change to polar coordinates
\begin{align*}
\xi &=r\begin{pmatrix}
\sin(\varphi)\sin(\alpha)\\
\cos(\varphi)\sin(\alpha)\\
\cos(\alpha)\end{pmatrix}. 
\end{align*}
In particular, that means $\theta =\eta_3/|\eta|=\cos(\alpha)$. In polar coordinates, $\Gamma_n$ writes 
\begin{align*}
\tilde \Gamma_n &=\{(r,\varphi ,\alpha) : \ 1-\lambda^2\le  \cos^2(\alpha) (4-3\gamma_n)\le 1+\lambda^2 \}.
\end{align*}
Hence
\begin{align*}
\left|\Gamma_n\cap\{2^{-1}\le|\eta|\le2\}\right|
&=\int_{2^{n+1}}^{2^{n-1}}\int_0^{2\pi}\int_0^\pi
\mathbf{1}_{\tilde\Gamma_n}(r,\varphi,\alpha)
r^2\sin(\alpha)\,d\alpha\,d\varphi\,dr\\
&=2\pi2^{-3n}\int_{2^{n-1}}^{2^{n+1}} r^2
\int_0^\pi
\mathbf{1}_{\{1-\lambda^2\le
\cos^2(\alpha)(4-3\gamma_n(r))\le1+\lambda^2\}}
\sin(\alpha)\,d\alpha\,dr.
\end{align*}

For each fixed $|\eta|\in[2^{-1},2]$, we distinguish two cases.
If $4-3\gamma(r)\le1/2$, then
\begin{equation*}
 \cos^2(\alpha)\bigl(4-3\gamma(r)\bigr)
\le\frac12<1-\lambda^2   
\end{equation*}
for all $\alpha\in[0,\pi]$. Therefore, the inner integral
vanishes.

If $4-3\gamma(r)>1/2$, the change of variables
$\theta=\cos(\alpha)$ yields
\begin{align*}
&\int_0^\pi
\mathbf{1}_{\{1-\lambda^2\le
\cos^2(\alpha)(4-3\gamma(r))\le1+\lambda^2\}}
\sin(\alpha)\,d\alpha\\
&\quad=\int_{-1}^{1}
\mathbf{1}_{\{1-\lambda^2\le
\theta^2(4-3\gamma(r))\le1+\lambda^2\}}
\,d\theta\\
&\quad\le
\frac{2}{\sqrt{4-3\gamma_n(r)}}
\left(\sqrt{1+\lambda^2}-\sqrt{1-\lambda^2}\right)\\
&\quad\le
2\sqrt{2}\,
\frac{2\lambda^2}
{\sqrt{1+\lambda^2}+\sqrt{1-\lambda^2}}
\lesssim\lambda^2.
\end{align*}
Here, 
\begin{equation*}
\sqrt{\frac{1-\lambda^2}{4-3\gamma(r)}}
\le|\theta|\le
\sqrt{\frac{1+\lambda^2}{4-3\gamma(r)}}.    
\end{equation*}

Intersecting these intervals with $[-1,1]$ can only decrease
their total length.

Combining the two cases, we obtain
\begin{equation*}
\left|\Gamma_n\cap\{2^{-1}\le|\eta|\le2\}\right|
\lesssim
\lambda^22^{-3n}\int_{2^{-1}}^2 r^2\,dr
\lesssim\lambda^2.    
\end{equation*}
\end{proof}
\begin{remark}
On $\chi(\theta)\neq0$, Lemma~\ref{lem:assd}
uniformly estimates, after dyadic rescaling, the volume
of the region where
\begin{equation*}
|1-\theta^2(4-3\gamma(r))|    
\end{equation*}
is small. This does not bound the number of stationary
points, which is controlled separately by
Assumption~\ref{ass:h}(B).

Since the full Hessian is singular on $\xi_3=0$,
the stationary set may contain the circle on
$\xi_3=0$, so an argument based on finitely many isolated
stationary points cannot be applied uniformly in the
spatial parameter.
We therefore obtain the $t^{-1}$ decay by applying
two-dimensional stationary phase on localized pieces.
\end{remark}
\section{Magnetic-Coriolis waves}\label{sec:mcw}

\begin{prop}[Precise version of Proposition 1]
For the initial data
$u_{0},b_{0}\in H^N(\mathbb R^3)$,
there exist
$\alpha_j,\beta_j$,
such that  the global in time solution $(u,b)$ of \eqref{eq:lin} is given by 
\begin{align*}
u(t,x)
&=\int_{\R^3}
\sum_{j=1}^4
\alpha_j(\xi)
e^{i(x\cdot\xi-\sigma_j(\xi)t)}
\,d\xi,\\
b(t,x)
&=\int_{\R^3}
\sum_{j=1}^4
\beta_j(\xi)
e^{i(x\cdot\xi-\sigma_j(\xi)t)}
\,d\xi,
\end{align*}
where $\sigma_j$ are given by \eqref{sigmaj}.
For $\xi\ne0$, choose  orthonormal vectors
$e_\pm(\xi)\in\xi^\perp\subset\mathbb C^3$ satisfying
\begin{equation*}
\Omega\left(
\mathrm{Id}-\frac{\xi\otimes\xi}{|\xi|^2}
\right)(e_3\wedge e_\pm)
=
\pm i\Omega\frac{\xi_3}{|\xi|}e_\pm.    
\end{equation*}
Then $\alpha_j,\beta_j$
are chosen as follows:
\begin{align*}
\alpha_1
&=
\frac{\sigma_1}{\sigma_1-\sigma_2}
\langle\widehat u_{0},e_+\rangle e_+
-
\frac{\xi\cdot B_0}{\sigma_1-\sigma_2}
\langle\widehat b_0,e_+\rangle e_+,\\
\alpha_2
&=
-\frac{\sigma_2}{\sigma_1-\sigma_2}
\langle\widehat u_{0},e_+\rangle e_+
+
\frac{\xi\cdot B_0}{\sigma_1-\sigma_2}
\langle\widehat b_0,e_+\rangle e_+,\\
\alpha_3
&=
\frac{\sigma_3}{\sigma_3-\sigma_4}
\langle\widehat u_{0},e_-\rangle e_-
-
\frac{\xi\cdot B_0}{\sigma_3-\sigma_4}
\langle\widehat b_0,e_-\rangle e_-,\\
\alpha_4
&=
-\frac{\sigma_4}{\sigma_3-\sigma_4}
\langle\widehat u_{0},e_-\rangle e_-
+
\frac{\xi\cdot B_0}{\sigma_3-\sigma_4}
\langle\widehat b_0,e_-\rangle e_-,
\end{align*}
and
\begin{align*}
\beta_1
&=
-\frac{\xi\cdot B_0}{\sigma_1-\sigma_2}
\langle\widehat u_{in},e_+\rangle e_+
-
\frac{\sigma_2}{\sigma_1-\sigma_2}
\langle\widehat b_0,e_+\rangle e_+,\\
\beta_2
&=
\frac{\xi\cdot B_0}{\sigma_1-\sigma_2}
\langle\widehat u_{0},e_+\rangle e_+
+
\frac{\sigma_1}{\sigma_1-\sigma_2}
\langle\widehat b_0,e_+\rangle e_+,\\
\beta_3
&=
-\frac{\xi\cdot B_0}{\sigma_3-\sigma_4}
\langle\widehat u_{0},e_-\rangle e_-
-
\frac{\sigma_4}{\sigma_3-\sigma_4}
\langle\widehat b_0,e_-\rangle e_-,\\
\beta_4
&=
\frac{\xi\cdot B_0}{\sigma_3-\sigma_4}
\langle\widehat u_{0},e_-\rangle e_-
+
\frac{\sigma_3}{\sigma_3-\sigma_4}
\langle\widehat b_0,e_-\rangle e_-.
\end{align*}

These formulas define the coefficients almost everywhere
in $\xi$. On the null set where the formulas are undefined,
including $\xi=0$, we set all coefficients equal to zero.
\end{prop}

\begin{proof}
First we derive the four frequencies by taking
the Fourier transform of the projected linear system.
Then we construct the coefficients and verify
the resulting solution formula.

\medskip
Applying the Leray projection for \eqref{eq:lin} and taking the
Fourier transform in space, we obtain

\begin{equation*}
 \begin{cases}
\partial_t\widehat u
+2\Omega P(\xi)(e_3\wedge\widehat u)
=i(B_0\cdot\xi)\widehat b,\\
\partial_t\widehat b
=i(B_0\cdot\xi)\widehat u,
\end{cases}
\qquad \text{where}\qquad
P(\xi)=\mathrm{Id}-\frac{\xi\otimes\xi}{|\xi|^2}.   
\end{equation*}
Since the $(u,b)$ are divergence-free, we write
\begin{equation*}
\widehat u=U_+e_++U_-e_-,
\qquad
\widehat b=B_+e_++B_-e_-.    
\end{equation*}
Then we turn to the system to
\begin{equation*}
    \begin{cases}
\partial_tU_++2i\Omega\frac{\xi_3}{|\xi|}U_+=iB_0\cdot\xi B_+,\\
\partial_tB_+=iB_0\cdot\xi U_+,
\end{cases}
\qquad
\begin{cases}
\partial_tU_--2i\Omega\frac{\xi_3}{|\xi|}U_-=iB_0\cdot\xi B_-,\\
\partial_tB_-=iB_0\cdot\xi U_-.
\end{cases}
\end{equation*}

For the first system, find a nonzero solution of the form
\begin{equation*}
 U_+(t)=A e^{-it\sigma},
\qquad
B_+(t)=B e^{-it\sigma}.   
\end{equation*}
Substitution gives
\begin{equation*}
    \begin{pmatrix}
2\Omega\frac{\xi_3}{|\xi|}-\sigma&-B_0\cdot\xi\\
-B_0\cdot\xi&-\sigma
\end{pmatrix}
\begin{pmatrix}A\\B\end{pmatrix}
=0.
\end{equation*}
A nonzero solution exists precisely when
\begin{equation*}
 (2\Omega\frac{\xi_3}{|\xi|}-\sigma)(-\sigma)-(B_0\cdot\xi)^2=0,   
\end{equation*}
or equivalently,
\begin{equation*}
 \sigma^2-2\Omega\frac{\xi_3}{|\xi|}\sigma-(B_0\cdot\xi)^2=0.   
\end{equation*}
Its roots are
\begin{equation*}
 \sigma_1=\Omega\frac{\xi_3}{|\xi|}+\sqrt{(\Omega\frac{\xi_3}{|\xi|})^2+(B_0\cdot\xi)^2},
\qquad
\sigma_2=\Omega\frac{\xi_3}{|\xi|}-\sqrt{(\Omega\frac{\xi_3}{|\xi|})^2+(B_0\cdot\xi)^2}.   
\end{equation*}
For the second system, the same calculation gives
\begin{equation}\label{sigma1234}
    \sigma_3=-\Omega\frac{\xi_3}{|\xi|}-\sqrt{(\Omega\frac{\xi_3}{|\xi|})^2+(B_0\cdot\xi)^2}=-\sigma_1,
\qquad
\sigma_4=-\Omega\frac{\xi_3}{|\xi|}+\sqrt{(\Omega\frac{\xi_3}{|\xi|})^2+(B_0\cdot\xi)^2}=-\sigma_2.
\end{equation}
These are the four frequencies appearing
in the solution formula.

\medskip
We first consider the $e_+$ and assume
$\sigma_1\ne\sigma_2$. Set
\begin{equation*}
\widehat u_{0,+}
=\langle\widehat u_0,e_+\rangle e_+,
\qquad
\widehat b_{0,+}
=\langle\widehat b_0,e_+\rangle e_+.    
\end{equation*}
Then
\begin{equation*}
\begin{aligned}
U_+(t)e_+
&=\alpha_1e^{-it\sigma_1}
+\alpha_2e^{-it\sigma_2},\\
B_+(t)e_+
&=\beta_1e^{-it\sigma_1}
+\beta_2e^{-it\sigma_2}.
\end{aligned}    
\end{equation*}
with $\alpha_1,\alpha_2,\beta_1,\beta_2$ to be determined.

Using the first system, gives
\begin{equation*}
\begin{cases}
\alpha_1+\alpha_2=\widehat u_{0,+},\\
\sigma_1\alpha_1+\sigma_2\alpha_2
=
2\Omega\frac{\xi_3}{|\xi|}
\widehat u_{0,+}
-(B_0\cdot\xi)\widehat b_{0,+},
\end{cases}    
\end{equation*}
and
\begin{equation*}
 \begin{cases}
\beta_1+\beta_2=\widehat b_{0,+},\\
\sigma_1\beta_1+\sigma_2\beta_2
=-(B_0\cdot\xi)\widehat u_{0,+}.
\end{cases}   
\end{equation*}
Then we can solve by
\begin{equation*}
 \alpha_1
=
\frac{
\sigma_1\widehat u_{0,+}
-(B_0\cdot\xi)\widehat b_{0,+}
}{\sigma_1-\sigma_2},
\qquad
\beta_1
=
\frac{
-(B_0\cdot\xi)\widehat u_{0,+}
-\sigma_2\widehat b_{0,+}
}{\sigma_1-\sigma_2}.   
\end{equation*}
The initial conditions then give
\begin{equation*}
\alpha_2
=
\frac{
-\sigma_2\widehat u_{0,+}
+(B_0\cdot\xi)\widehat b_{0,+}
}{\sigma_1-\sigma_2},
\qquad
\beta_2
=
\frac{
(B_0\cdot\xi)\widehat u_{0,+}
+\sigma_1\widehat b_{0,+}
}{\sigma_1-\sigma_2}.    
\end{equation*}
The coefficients $\alpha_3,\alpha_4,\beta_3,\beta_4$
are obtained in the same way from the $e_-$ system.

Then we obtain
\begin{equation*}
|\alpha_j(\xi)|+|\beta_j(\xi)|
\le C\bigl(
|\widehat u_0(\xi)|+|\widehat b_0(\xi)|
\bigr),
\qquad j=1,\ldots,4.    
\end{equation*}
 Multiplying by
$(1+|\xi|^2)^{N/2}$ and taking $L^2$ norms, we obtain
\begin{equation}\label{eq:coefficients-HN}
\sum_{j=1}^4
\left(
\|\mathcal F^{-1}\alpha_j\|_{H^N}
+
\|\mathcal F^{-1}\beta_j\|_{H^N}
\right)
\le
C\left(\|u_0\|_{H^N}+\|b_0\|_{H^N}\right).
\end{equation}

Define
\begin{equation*}
\widehat u(t,\xi)
=\sum_{j=1}^4\alpha_j(\xi)e^{-it\sigma_j(\xi)},
\qquad
\widehat b(t,\xi)
=\sum_{j=1}^4\beta_j(\xi)e^{-it\sigma_j(\xi)}.    
\end{equation*}
Estimate \eqref{eq:coefficients-HN}
and the dominated convergence theorem therefore give \begin{equation*}
 (u,b)\in C(\mathbb R;H^N\times H^N).   
\end{equation*}
The coefficient formulas and the dispersion relations
give, 
\begin{equation*}
\begin{aligned}
-i\sigma_j\alpha_j
+2\Omega P(\xi)(e_3\wedge\alpha_j)
-i(B_0\cdot\xi)\beta_j&=0,\\
\sigma_j\beta_j+(B_0\cdot\xi)\alpha_j&=0.
\end{aligned}    
\end{equation*}
Then by constructing
\begin{equation*}
\sum_{j=1}^4\alpha_j=\widehat u_0,
\qquad
\sum_{j=1}^4\beta_j=\widehat b_0.    
\end{equation*}
Finally, the projected system satisfies
\begin{equation*}
\|u(t)\|_{H^N}^2+\|b(t)\|_{H^N}^2=
\|u_0\|_{H^N}^2+\|b_0\|_{H^N}^2.   
\end{equation*}
Recall the pressure gradient through
\begin{equation*}
\nabla\pi
=
-2\Omega(\mathrm{Id}-\mathbb P)(e_3\wedge u)    
\end{equation*}
shows that $(u,b)$ solves \eqref{eq:lin}
with the prescribed initial data.
\end{proof}

\begin{lemma}\label{lem:ms}
Let $\alpha,\Omega\ne 0$ and define
\begin{equation*}
 h_i(r)
:=
\frac{
1+(-1)^i\left\langle\frac{\alpha}{\Omega}r\right\rangle
}{
(\frac{\alpha}{\Omega})r
},
\qquad i\in\{0,1\},\quad r>0,
\qquad
\langle s\rangle:=\sqrt{1+s^2}.   
\end{equation*}

Then the profiles $\alpha h_i(r)$ satisfy
Assumption~\ref{ass:h} with
\begin{equation*}
\zeta_0=\zeta_\infty=\zeta
:=
\frac{\alpha^2|\Omega|}{\Omega^2+\alpha^2},    
\end{equation*}
and
\begin{equation*}
m_0=m_\infty=2\quad\text{if }i=0,
\qquad
m_0=0,\quad m_\infty=2\quad\text{if }i=1.    
\end{equation*}
The constants are uniform in $\alpha$ and $\Omega$.
One may take $a_0=\frac34$ and $A_0=3$.
In condition (B), the stationary point counts are
at most two for $i=0$ and at most four for $i=1$.
\end{lemma}

The proof of this lemma can be found in the appendix. With this lemma, we can infer the decay of the magnetic Coriolis waves using Theorem \ref{prop:strich}:

\begin{cor}[Precise version of Corollary \ref{linearcase}]\label{cor:decay}
For $N\ge 4 $, let the divergence-free initial data satisfy
\begin{equation*}
u_{0},b_{0}\in
H^N(\R^3)
\cap\dot B^1_{1,1}(\R^3)
\cap\dot B^3_{1,1}(\R^3).\end{equation*}
then for the solution of the linearized equations \eqref{eq:lin}, it holds 
\begin{align*}
\Vert (u,b)(t)\Vert_{L^\infty}&\le C\frac 1 { t  \zeta +1} \Vert (u_{0},b_{0})\Vert_{\dot B_{1,1}^{1}\cap \dot B_{1,1}^{3}\cap H^N}. \end{align*}
In particular, let \(\mathcal L(t)\) denote the solution operator associated with System \eqref{eq:lin}, namely
\begin{align*}
\mathcal L(t)(u_{0},b_{0})=(u(t),b(t)).    
\end{align*}
Then
\begin{align*}
\|\mathcal L(t)\|_{X\to L^\infty}
\lesssim
\frac 1 { t  \zeta +1},  
\end{align*}
where
$
X=\left\{
(u,b):
u,b\in
\dot B_{1,1}^{1}(\mathbb R^3;\mathbb R^3)
\cap
\dot B_{1,1}^{3}(\mathbb R^3;\mathbb R^3)\cap H^N(\R^3;\R^3),
\quad
\nabla\cdot u=\nabla\cdot b=0
\right\}.$
\end{cor}
\begin{proof}
On $\{\Omega\xi_3<0\}$, we interchange
$\sigma_1$ by $\sigma_2$ and $\sigma_3$ by $\sigma_4$
in the phases and their corresponding coefficients,
and retain the same notation. Recall that
\begin{equation*}
 \begin{aligned}
\sigma_1(\xi)
&=\alpha\xi_3h_0(|\xi|),
& \sigma_3&=-\sigma_1,\\
\sigma_2(\xi)
&=\alpha\xi_3h_1(|\xi|),
& \sigma_4&=-\sigma_2.
\end{aligned}   
\end{equation*}
By Lemma~\ref{lem:ms} and
$(\alpha h_i)'=\alpha h_i'$, the profiles
$\pm\alpha h_i$ satisfy Assumption~\ref{ass:h}.
Since $\zeta\le\min\{\frac{\alpha^2}{\Omega},\Omega\},$
both profiles satisfy condition~(A) with
$\zeta_0=\zeta_\infty=\zeta$, and
\begin{equation*}
    (m_0,m_\infty)=
\begin{cases}
(2,2),&i=0,\\
(0,2),&i=1.
\end{cases}
\end{equation*}

Applying Theorem~\ref{prop:strich},
we obtain, for $j=1,3$ and $t>0$,
\begin{align*}
&\left\|\mathcal F^{-1}
[e^{-it\sigma_j}\alpha_j]\right\|_{L^\infty}
+
\left\|\mathcal F^{-1}
[e^{-it\sigma_j}\beta_j]\right\|_{L^\infty}\\
\lesssim&
\frac1{\zeta t}
\left(
\|\mathcal F^{-1}\alpha_j\|_{\dot B^3_{1,1}}
+
\|\mathcal F^{-1}\beta_j\|_{\dot B^3_{1,1}}
\right),\\
\lesssim&
\|(u_0,b_0)\|_{\dot B^s_{1,1}},
\qquad \quad s\in\{1,3\}.
\end{align*}
For $j=2,4$, thanks to Theorem~\ref{prop:strich} again
\begin{align*}
&\left\|\mathcal F^{-1}
[e^{-it\sigma_j}\alpha_j]\right\|_{L^\infty}
+
\left\|\mathcal F^{-1}
[e^{-it\sigma_j}\beta_j]\right\|_{L^\infty}\\
&\qquad\lesssim
\frac1{\zeta t}
\sum_{s\in\{1,3\}}
\left(
\|\mathcal F^{-1}\alpha_j\|_{\dot B^s_{1,1}}
+
\|\mathcal F^{-1}\beta_j\|_{\dot B^s_{1,1}}
\right).
\end{align*}
Summing over $j$ and using the coefficient bounds
\begin{equation}\label{eq:linear-large-time}
\|(u,b)(t)\|_{L^\infty}
\lesssim
\frac{1}{\zeta t}
\left(
\|(u_0,b_0)\|_{\dot B^1_{1,1}}
+
\|(u_0,b_0)\|_{\dot B^3_{1,1}}
\right),
\qquad t>0.
\end{equation}
we conclude that
\begin{equation}\label{eq:linear-large-time}
\|(u,b)(t)\|_{L^\infty}
\lesssim
\frac1{\zeta t}
\|(u_0,b_0)\|_{\dot B^1_{1,1}\cap\dot B^3_{1,1}},
\qquad t>0.
\end{equation}

On the other hand, Sobolev embedding and energy
conservation give
\begin{equation*}
\begin{aligned}
\|u(t)\|_{L^\infty}+\|b(t)\|_{L^\infty}
&\lesssim
\left(\|u(t)\|_{H^N}^2+\|b(t)\|_{H^N}^2\right)^{1/2}\\
&=
\left(\|u_0\|_{H^N}^2+\|b_0\|_{H^N}^2\right)^{1/2},
\end{aligned}    
\end{equation*}
where $N\ge4$.
Combining this bound with \eqref{eq:linear-large-time},
we obtain, for $t>0$,
\begin{equation*}
 \begin{aligned}
\|(u,b)(t)\|_{L^\infty}
&\lesssim
\min\left\{1,\frac1{\zeta t}\right\}
\|(u_0,b_0)\|_
{\dot B^1_{1,1}\cap\dot B^3_{1,1}\cap H^N}\\
&\lesssim
\frac1{1+\zeta t}
\|(u_0,b_0)\|_
{\dot B^1_{1,1}\cap\dot B^3_{1,1}\cap H^N}.
\end{aligned}   
\end{equation*}
The estimate at $t=0$ follows from Sobolev embedding.
\end{proof}

\subsection{Extension of time of existence}
In this section, we apply the dispersive estimates to the nonlinear equations. For a rapidly rotating MHD system with vanishing viscosity,
\cite{CharveYao} combined Strichartz estimates
with energy methods to obtain global strong solutions
and explicit convergence rates.In the present work, we combine the linear decay estimate
in Corollary \ref{cor:decay} with Sobolev energy estimates
to obtain the explicit lifespan bound stated in
Theorem \ref{thm:extend}.

\begin{prop}\label{bkm}
Let $u, b \in C([0,T],H^N(\mathbb{R}^3))$ with $N > \tfrac{5}{2}$ then there exists a constant $C>0$ such that 
\begin{equation*}
\Vert u, b\Vert_{H^N}(t)\le \Vert u_0, b_0\Vert_{H^N}\exp\left(C\int_0^t\Vert (\nabla(u,b)(\tau)\Vert_{L^\infty}d\tau\right)
\end{equation*}
for all $t\le T$.
\end{prop}
Since this is a standard Beale-Kato-Majda criterion, we omit the proof and refer to \cite{caflisch1997remarks,elgindi2015sharp}. Using this, we can prove the extended time of existence:
\begin{proof}[Proof of Theorem~\ref{thm:extend}]
Let $T_{\max}>0$ be the maximal existence time.
By Proposition~\ref{bkm},
\begin{equation}\label{eq:energy-extension}
\|(u,b)(t)\|_{H^N}
\le
\varepsilon
\exp\left(
C\int_0^t
\|\nabla(u,b)(s)\|_{L^\infty}\,ds
\right),
\qquad 0\le t<T_{\max}.
\end{equation}

Since
\begin{equation*}
\|(u_0,b_0)\|_{H^N}=\varepsilon<2\varepsilon    
\end{equation*}
and $(u,b)$ is continuous in time with values in $H^N$,
there exists $T_0\in(0,T_{\max})$ such that
\begin{equation*}
    \sup_{0\le t\le T_0}\|(u,b)(t)\|_{H^N}<2\varepsilon.
\end{equation*}
Thus the bootstrap bound holds.

Now let $0<T<T_{\max}$ such that
\begin{equation}\label{eq:bootstrap-extension}
\sup_{0\le t\le T}\|(u,b)(t)\|_{H^N}
\le2\varepsilon.
\end{equation}
We establish estimates on this interval and then determine
how large $T$ can be while strictly improving this bound.

By Duhamel's formula,
\begin{equation*}
\begin{aligned}
\begin{pmatrix}u\\ b\end{pmatrix}(t)
&=
\mathcal L(t)\begin{pmatrix}u_0\\ b_0\end{pmatrix}\\
&\quad+
\int_0^t\mathcal L(t-s)
\begin{pmatrix}
\mathbb P\bigl((b\cdot\nabla)b-(u\cdot\nabla)u\bigr)\\
(b\cdot\nabla)u-(u\cdot\nabla)b
\end{pmatrix}(s)\,ds.
\end{aligned}
\end{equation*}
The linear evolution commutes with spatial derivatives.
Combining the dispersive estimate with Sobolev embedding
and conservation of the linear $H^4$ energy gives
\begin{equation*}
 \|\nabla\mathcal L(t)f\|_{L^\infty}
\lesssim
\frac1{1+\zeta t}
\|f\|_{\dot B^2_{1,1}\cap\dot B^4_{1,1}\cap H^4},
\qquad t\ge0,   
\end{equation*}
for divergence-free $f$.
Using also the boundedness of the Leray projection
on these spaces, we obtain
\begin{equation*}
\begin{aligned}
\|\nabla(u,b)(t)\|_{L^\infty}
&\le
\frac{CM_0}{1+\zeta t}\\
&\quad+
C\sum_{a^1,a^2\in\{u,b\}}
\int_0^t
\frac{
\|(a^1\cdot\nabla)a^2(s)\|_
{\dot B^2_{1,1}\cap\dot B^4_{1,1}\cap H^4}
}{
1+\zeta(t-s)
}\,ds,
\end{aligned}    
\end{equation*}
where
\begin{equation*}
M_0:=
\|(u_0,b_0)\|_
{\dot B^2_{1,1}\cap\dot B^4_{1,1}\cap H^4}.
\end{equation*}
To estimate the nonlinear terms, we use Proposition\ref{embedding}
\begin{equation*}
W^{5,1}(\mathbb R^3)
\hookrightarrow
\dot B^2_{1,1}(\mathbb R^3)
\cap\dot B^4_{1,1}(\mathbb R^3).    
\end{equation*}
The Hölder inequality imply
\begin{equation*}
\|(a^1\cdot\nabla)a^2\|_{W^{5,1}}
\lesssim
\|a^1\|_{H^5}\|a^2\|_{H^6}.
\end{equation*}
Moreover, Proposition\ref{sobolev-product} of $H^4(\mathbb R^3)$
gives
\begin{equation*}
    \|(a^1\cdot\nabla)a^2\|_{H^4}
\lesssim
\|a^1\|_{H^4}\|a^2\|_{H^5}.
\end{equation*}
Since $N>6$, it follows that
\begin{equation*}
    \|(a^1\cdot\nabla)a^2\|_
{\dot B^2_{1,1}\cap\dot B^4_{1,1}\cap H^4}
\lesssim
\|a^1\|_{H^N}\|a^2\|_{H^N}.
\end{equation*}
Hence \eqref{eq:bootstrap-extension} yields
\begin{equation*}
\begin{aligned}
\|\nabla(u,b)(t)\|_{L^\infty}
&\le
\frac{CM_0}{1+\zeta t}
+
C\varepsilon^2
\int_0^t\frac{ds}{1+\zeta(t-s)}\\
&=
\frac{CM_0}{1+\zeta t}
+
\frac{C\varepsilon^2}{\zeta}
\log(1+\zeta t).
\end{aligned}    
\end{equation*}
Integrating over $[0,T]$, we obtain
\begin{equation}\label{eq:integrated-gradient}
\int_0^T\|\nabla(u,b)(t)\|_{L^\infty}\,dt
\le
C\left(
\frac{M_0}{\zeta}
+
\frac{\varepsilon^2T}{\zeta}
\right)\log(1+\zeta T).
\end{equation}

To improve the bootstrap bound using
\eqref{eq:energy-extension}, we need the right-hand side
of \eqref{eq:integrated-gradient} to be sufficiently small.
The nonlinear contribution is controlled by
\begin{equation*}
 \frac{\varepsilon^2T}{\zeta}\log(1+\zeta T).   
\end{equation*}

This leads us to set
\begin{equation*}
T_*:=\frac{c\zeta}{\varepsilon^2\log(2+\zeta^2\varepsilon^{-2})},
\end{equation*}
where $0<c\le1$ will be chosen sufficiently small.

The assumptions on the initial data and the embedding give
\begin{equation*}
    \begin{aligned}
M_0
&\le
C\left(
\|(u_0,b_0)\|_{W^{N-1,1}}
+
\|(u_0,b_0)\|_{H^N}
\right)\\
&\le
C\max\left\{\varepsilon,\frac{c_0\zeta}{\log(2+\zeta^2\varepsilon^{-2})}\right\}.
\end{aligned}
\end{equation*}
Also,
\begin{equation*}
    \frac{\varepsilon}{\zeta}<c_0.
\end{equation*}
Choose $c_0\le e^{-1}$. Then $\log(2+\zeta^2\varepsilon^{-2})>1$, and for $T\le T_*$,
\begin{equation*}
\log(1+\zeta T)
\le
\log\left(
1+\frac{c\zeta^2}{\varepsilon^2\log(2+\zeta^2\varepsilon^{-2})}
\right)
\le \log(2+\zeta^2\varepsilon^{-2}).
\end{equation*}
Consequently,
\begin{equation*}
    \frac{\varepsilon^2T}{\zeta}\log(1+\zeta T)
\le c,
\end{equation*}
and
\begin{equation*}
    \frac{M_0}{\zeta}\log(1+\zeta T)
\le
C\max\left\{\frac{\varepsilon}{\zeta}\log(2+\zeta^2\varepsilon^{-2}),c_0\right\}.
\end{equation*}

Since $\frac{\varepsilon}{\zeta}<c_0\le e^{-1}$,
\begin{equation*}
    \log(2+\zeta^2\varepsilon^{-2})\le
\log3+2\log\left(\frac{\zeta}{\varepsilon}\right).
\end{equation*}
We know that on
$(0,e^{-1}]$ therefore implies
\begin{equation*}
\frac{\varepsilon}{\zeta}\log(2+\zeta^2\varepsilon^{-2})
\le
c_0\log3+2c_0\log(\frac1{c_0}).
\end{equation*}
Substituting these bounds into
\eqref{eq:integrated-gradient}, we obtain
\begin{equation}\label{eq:small-gradient-integral}
\int_0^T\|\nabla(u,b)(t)\|_{L^\infty}\,dt
\le
C\left[c_0(1+\log(\frac1{c_0}))+c\right]
\end{equation}
whenever $T<T_{\max}$, $T\le T_*$,
and \eqref{eq:bootstrap-extension} holds.

By \eqref{eq:energy-extension},
\begin{equation*}
\sup_{0\le t\le T}\|(u,b)(t)\|_{H^N}
\le
\varepsilon
\exp\left(
C\left[c_0(1+\log(\frac1{c_0}))+c\right]
\right).
\end{equation*}
Choosing $c_0$ and then $c$ sufficiently small gives
\begin{equation}\label{eq:bootstrap-improvement}
\sup_{0\le t\le T}\|(u,b)(t)\|_{H^N}
\le\frac32\varepsilon.
\end{equation}

If the norm reached $2\varepsilon$ at some time
$t<\min\{T_{\max},T_*\}$, time continuity and the
initial strict bound would give a first such time $t_1$.
On $[0,t_1]$, the bootstrap bound would hold, so
\eqref{eq:bootstrap-improvement} would imply
\begin{equation*}
    2\varepsilon
=\|(u,b)(t_1)\|_{H^N}
\le\frac32\varepsilon,
\end{equation*}
a contradiction.
Thus the bootstrap bound holds throughout
$[0,\min\{T_{\max},T_*\})$.

If $T_{\max}\le T_*$, letting $T$ tends to $T_{\max}$
in \eqref{eq:small-gradient-integral} gives
\begin{equation*}
    \int_0^{T_{\max}}
\|\nabla(u,b)(t)\|_{L^\infty}\,dt<\infty.
\end{equation*}
Proposition~\ref{bkm} then extends the solution beyond
$T_{\max}$, contradicting maximality.
Hence $T_{\max}>T_*$.
By time continuity,
\begin{equation*}
    \sup_{0\le t\le T_*}\|(u,b)(t)\|_{H^N}
\le\frac32\varepsilon<2\varepsilon.
\end{equation*}
Taking $T=T_*$ proves the theorem.
\end{proof}

\appendix
\section{Technical Lemmas}
\textbf{Proof of Lemma \ref{lem:ms}:}
Let
\begin{equation*}
h_i(r):=\frac{1+(-1)^i\left\langle \frac{\alpha}{\Omega}r\right\rangle}{(\frac{\alpha}{\Omega})r},
\qquad i\in\{0,1\},
\qquad r\neq 0,   
\end{equation*}
where
\begin{equation*}
 \langle s\rangle:=\sqrt{1+s^2}.  
\end{equation*}

Let $\chi$ be the angular cutoff fixed in Assumption~4.
Since $0\notin\supp\chi$, there exists $\delta>0$ such that
\begin{equation*}
\supp\chi\subset
\left\{s\in[-1,1]:|s|\ge\delta\right\}.    
\end{equation*}
Thus, $\chi(\theta)$ localizes away from the equatorial
plane $\theta=\xi_3/|\xi|=0$, where the stationary set
may fail to be discrete for certain values of the
spatial parameter.
Then, for each $i\in\{0,1\}$, there exists a uniform constant $c>0$
such that the following hold:
\begin{itemize}
\item If $i=0$, then condition \textnormal{(A)} in Assumption \ref{ass:h} holds with
\begin{equation*}
m_0=m_\infty=2,
\qquad
\zeta_0=\zeta_\infty=\frac{\Omega}{\alpha}.    
\end{equation*}
\item If $i=1$, then condition \textnormal{(A)} in Assumption \ref{ass:h} holds with
\begin{equation*}
  m_0=0,\qquad m_\infty=2,
    \qquad
    \zeta_0=\frac{\alpha}{\Omega},
    \qquad
\zeta_\infty=\frac{\Omega}{\alpha}.   
\end{equation*}

\item Condition \textnormal{(B)} in Assumption \ref{ass:h} holds, namely
\begin{equation*}
   \sup_{x,n}
    \#\Bigl\{
    \xi:\ 2^n\le |\xi|\le 2^{n+1},\ \chi(\xi_3/|\xi|)\neq 0,\ 
    \nabla_\xi \bigl(x\cdot \xi-\xi_3 h_i(|\xi|)\bigr)=0
    \Bigr\}
    \le c^{-1}.  
\end{equation*}
In fact, one may take the right-hand side to be $2$ (when $i=0$), and to be $4$ (when $i=1$).
\end{itemize}
With out loss of generality, we set that  $a:=\frac{\alpha}{\Omega}>0,$ then
\begin{equation*}
h_i(r)=\frac{1+(-1)^i\langle ar\rangle}{ar}.    
\end{equation*}
\subsection{}
For $i=0$, we compute
\begin{equation*}
\begin{aligned}
\vert h_0'(r)\vert =\frac{1}{ar^2}\Bigl(1+\frac1{\langle ar\rangle}\Bigr).
    \end{aligned}
\end{equation*}

And we have
\begin{equation*}
\frac{1}{ar^2}\le |h_0'(r)|\le \frac{2}{ar^2}.    
\end{equation*}
It follows that
\begin{equation*}
|h_0'(r)|\ge c\,\frac{\Omega}{\alpha}\, r^{-2}.   
\end{equation*}
Thus condition \textnormal{(A)} holds for $h_0$ with
$
m_0=m_\infty=2$ and $\zeta_0=\zeta_\infty=\frac{\Omega}{\alpha}$.

For $i=1$, here we have 
\begin{equation*}
|h_1'(r)|=\frac{a}{\langle ar\rangle(1+\langle ar\rangle)}.   \end{equation*}
We now split the regions $|r|\le a^{-1}$ and $|r|\ge a^{-1}$.

If $|r|\le a^{-1}$, then $\langle ar\rangle\approx 1$, and therefore
\begin{equation*}
|h_1'(r)|\approx a=\frac{\alpha}{\Omega}.    
\end{equation*}
Thus
\begin{equation*}
|h_1'(r)|\ge c\,\frac{\alpha}{\Omega}
=
c\,\zeta_0\, r^{-m_0}
\qquad\text{with}\qquad
m_0=0,\ \ \zeta_0=\frac{\alpha}{\Omega}.    
\end{equation*}

If $|r|\ge a^{-1}$, then $\langle ar\rangle\approx a|r|$, and so
\begin{equation*}
|h_1'(r)|
\approx
\frac{a}{(a|r|)(a|r|)}=\frac{1}{a r^2}=
\frac{\Omega}{\alpha}\,r^{-2}.    
\end{equation*}
Hence
\begin{equation*}
 |h_1'(r)|
\ge
c\,\frac{\Omega}{\alpha}\,r^{-2}
=
c\,\zeta_\infty\, r^{-m_\infty}
\qquad\text{with}\qquad
m_\infty=2,\ \ \zeta_\infty=\frac{\Omega}{\alpha}.   
\end{equation*}

Combining the two regions, we conclude that
\begin{equation*}
 |h_1'(r)|
\ge
c\,\min\!\left(\frac{\alpha}{\Omega},\frac{\Omega}{\alpha}r^{-2}\right),   
\end{equation*}
which is exactly condition \textnormal{(A)} for $h_1$ with $m_0=0$, $ m_\infty=2$, 
$\zeta_0=\frac{\alpha}{\Omega}$, and 
$\zeta_\infty=\frac{\Omega}{\alpha}$.

\subsection{}
For any $x\in\R^3$ and we consider
\begin{equation*}
\Phi_x(\xi):=x\cdot \xi-\xi_3 h_i(|\xi|).    
\end{equation*}
Since
\begin{equation*}
\nabla(\xi_3 h_i(r))
=
h_i(r)e_3+\xi_3 h_i'(r)\frac{\xi}{r},    
\end{equation*}
and
\begin{equation*}
\nabla\varphi_x(\xi)=0    
\end{equation*}
is equivalent to
\begin{equation*}
x=h_i(r)e_3+\xi_3 h_i'(r)\frac{\xi}{r}.    
\end{equation*}
Using the notations$\xi=(\xi_h,\xi_3)$ and $x=(x_h,x_3)$, we obtain
\begin{equation*}
x_h=\xi_3 h_i'(r)\frac{\xi_h}{r},
\qquad
x_3=h_i(r)+\frac{\xi_3^2}{r}h_i'(r).
\end{equation*}
Since $|\xi_h|^2=r^2(1-\theta^2)$, this implies
\begin{equation}\label{cshzeta}
|x_h|^2=(r\theta h_i'(r))^2 (1-\theta^2),
\end{equation}
and
\begin{equation}\label{xhzeta}
x_3=h_i(r)+\theta^2 r h_i'(r).
\end{equation}
Then we obtain
\begin{equation}\label{fhzeta}
|x_h|^2
=
F_i(r):=
\bigl(x_3-h_i(r)\bigr)\bigl(r h_i'(r)-x_3+h_i(r)\bigr).
\end{equation}
Differentiating \eqref{fhzeta} and using \eqref{xhzeta}, we obtain
\begin{equation*}
F_i'(r)
=
r\,h_i'(r)^2
\Bigl(
\theta^2\Bigl(3+\frac{r h_i''(r)}{h_i'(r)}\Bigr)-1
\Bigr).    
\end{equation*}

For $i=0$, we compute that
\begin{align}\label{hfrac0}
\frac{r h_0''(r)}{h_0'(r)}
=
-2-\frac1{\langle ar\rangle}+\frac1{\langle ar\rangle^2}, 
\end{align}
hence
\begin{equation*}
\frac34<3+\frac{r h_0''(r)}{h_0'(r)}
=
1-\frac1{\langle ar\rangle}+\frac1{\langle ar\rangle^2}<1.    
\end{equation*}
Therefore
\begin{equation*}
\theta^2\Bigl(3+\frac{r h_0''(r)}{h_0'(r)}\Bigr)-1<0.   
\end{equation*}

For $i=1$, $F_1'$ does not have a fixed sign.
So we define
\begin{align*}
 z:=\frac{ar}{1+\langle ar\rangle}\in(0,1),   
\end{align*}
it's strictly increasing, and
\begin{align*}
h_1(r)=-z,
\qquad
h_1(r)+rh_1'(r)
=
-\frac{2z}{1+z^2}.    
\end{align*}
 \eqref{fhzeta} becomes
 \begin{align*}
|x_h|^2
=
-(x_3+z)
\left(
x_3+\frac{2z}{1+z^2}
\right).     
 \end{align*}
Moreover,
\begin{align*}
\frac{d}{dz}F_1(r(z))
=
-\frac{x_3(3+z^4)+4z}{(1+z^2)^2}. \end{align*}
And we know that the function $F_1(r(z))$ has at most one critical point.
Consequently, the equation $F_1(r)=|x_h|^2$
has at most two solutions $r>0$.
Once $r$ is fixed, \eqref{xhzeta} determines
\begin{equation*}
\theta^2=\frac{x_3-h_i(r)}{rh_i'(r)},   
\end{equation*}
and hence there are at most two possible values of $\theta$.
Since $\chi(\theta)\neq0$ implies $|\theta|\geq\delta$, we have
$\theta\neq0$. Then we obtain
\begin{align*}
  \xi_h=\frac{x_h}{\theta h_i'(r)}.  
\end{align*}

Therefore, for each fixed $x$ and each dyadic shell, there are at
most two stationary points for $h_0$ and at most four stationary
points for $h_1$. In particular,
\begin{align*}
 \sup_{x,n}
\#
\left\{
\xi:
2^n\leq|\xi|\leq2^{n+1},
\quad
\chi\left(\frac{\xi_3}{|\xi|}\right)\neq0,
\quad
\nabla_\xi
\left(
x\cdot\xi-\xi_3h_i(|\xi|)
\right)=0
\right\}
\leq4.   
\end{align*}
This proves condition \textnormal{(B)}.
\subsection{}
For $i=0$, 
\begin{equation*}
h_0'(r)= -\frac{1}{ar^2}\Big(1+\frac1{\langle ar\rangle}\Big),\,\,\,h_0''(r)
=
\frac{2}{ar^3}\Bigl(1+\frac1{\langle ar\rangle}\Bigr)
+\frac{a}{r\langle ar\rangle^3}.  \end{equation*}
Hence
\begin{equation*}
\frac{r h_0''(r)}{h_0'(r)}
=
-2-\frac{a^2r^2}{\langle ar\rangle^3(1+\langle ar\rangle^{-1})}.    
\end{equation*}
Since
\begin{equation*}
\frac{a^2r^2}{\langle ar\rangle^3(1+\langle ar\rangle^{-1})}
=
\frac{\langle ar\rangle^2-1}{\langle ar\rangle^2(\langle ar\rangle+1)}
=
\frac{\langle ar\rangle-1}{\langle ar\rangle^2}
=
\frac1{\langle ar\rangle}-\frac1{\langle ar\rangle^2},   
\end{equation*}
we have
\begin{equation*}
3+\frac{r h_0''(r)}{h_0'(r)}
=
1-\frac1{\langle ar\rangle}+\frac1{\langle ar\rangle^2}.    
\end{equation*}

Similarly, for $i=1$, we obtain
\begin{equation*}
 3+\frac{r h_1''(r)}{h_1'(r)}
=
1+\frac1{\langle ar\rangle}+\frac1{\langle ar\rangle^2}.   
\end{equation*}

Condition~\textnormal{(C)} holds with
$a_0=\frac34$ and $A_0=3$.

\subsection{}
Both profiles belong to $C^\infty((0,\infty))$.
Differentiating the preceding identities gives
\begin{equation*}
  \frac{r^2h_i'''(r)}{h_i'(r)}
=
6+\frac{6(-1)^i}{\langle ar\rangle}
-\frac6{\langle ar\rangle^2}
-\frac{3(-1)^i}{\langle ar\rangle^3}
+\frac3{\langle ar\rangle^4}.  
\end{equation*}
Consequently,
\begin{equation*}
r^2|h_i'''(r)|\le24|h_i'(r)|.
\end{equation*}

Thus (D) holds.

\subsection*{Data availability}
No data was used for the research described in the article.

\subsection*{Acknowledgments}

N. Knobel is supported by ERC/EPSRC Horizon Europe Guarantee EP/X020886/1. The authors declare that he has no conflict of interest.

The authors would like to thank Pierre Germain for useful discussions and thank Luis Miguel Rodrigues for recommending references.

\clearpage

\bibliographystyle{alpha} 
\bibliography{library}

\newcommand{\etalchar}[1]{$^{#1}$}
\begin{thebibliography}{HLR{\etalchar{+}}19}

\bibitem[Alf42]{alfven1942existence}
Hannes Alfv{\'e}n.
\newblock Existence of electromagnetic-hydrodynamic waves.
\newblock {\em Nature}, 150(3805):405--406, 1942.

\bibitem[BSS88]{bardos1988longtime}
Claude Bardos, Catherine Sulem, and Pierre-Louis Sulem.
\newblock Longtime dynamics of a conductive fluid in the presence of a strong
  magnetic field.
\newblock {\em Transactions of the American Mathematical Society},
  305(1):175--191, 1988.

\bibitem[CKS97]{caflisch1997remarks}
Russel~E Caflisch, Isaac Klapper, and Gregory Steele.
\newblock Remarks on singularities, dimension and energy dissipation for ideal
  hydrodynamics and mhd.
\newblock {\em Communications in Mathematical Physics}, 184(2):443--455, 1997.

\bibitem[CY26]{CharveYao}
Fr{\'e}d{\'e}ric Charve and Zhuoya Yao.
\newblock Convergence rates for a rotating {MHD} system with evanescent
  viscosity.
\newblock {\em arXiv preprint arXiv:2606.09354}, 2026.

\bibitem[CZ24]{chen2024sobolev}
Ting Chen and Ruizhao Zi.
\newblock On the {S}obolev stability threshold for shear flows near {C}ouette
  in {2D MHD} equations.
\newblock {\em Proceedings of the Royal Society of Edinburgh Section A:
  Mathematics}, pages 1--51, 2024.

\bibitem[Dav16]{davidson_2016}
Peter~A Davidson.
\newblock {\em Introduction to Magnetohydrodynamics}.
\newblock Cambridge Texts in Applied Mathematics. Cambridge University Press, 2
  edition, 2016.

\bibitem[DKZ24]{knobel2024nr}
Michele Dolce, Niklas Knobel, and Christian Zillinger.
\newblock {L}arge norm inflation of the current in the viscous, non-resistive
  magnetohydrodynamics equations.
\newblock {\em arXiv preprint arXiv:2410.22804}, 2024.

\bibitem[Dol24]{Dolce}
Michele Dolce.
\newblock Stability threshold of the {2D} {C}ouette flow in a homogeneous
  magnetic field using symmetric variables.
\newblock {\em Communications in Mathematical Physics}, 405(4):94, 2024.

\bibitem[EW15]{elgindi2015sharp}
Tarek~M Elgindi and Klaus Widmayer.
\newblock Sharp decay estimates for an anisotropic linear semigroup and
  applications to the surface quasi-geostrophic and inviscid boussinesq
  systems.
\newblock {\em SIAM Journal on Mathematical Analysis}, 47(6):4672--4684, 2015.

\bibitem[FMN19]{ferriz2019advances}
Antonio Ferriz-Mas and Manuel N{\'u}{\~n}ez.
\newblock {\em Advances in nonlinear dynamos}.
\newblock CRC Press, 2019.

\bibitem[GGV17]{gallagher2017wellposedness}
Isabelle Gallagher and David G{\'e}rard-Varet.
\newblock Wellposedness of linearized taylor equations in magnetohydrodynamics.
\newblock In {\em Shocks, Singularities and Oscillations in Nonlinear Optics
  and Fluid Mechanics}, pages 103--140. Springer, 2017.

\bibitem[GIP16]{guo2016global}
Yan Guo, Alexandru~D Ionescu, and Benoit Pausader.
\newblock Global solutions of the euler—maxwell two-fluid system in 3d.
\newblock {\em annals of Mathematics}, pages 377--498, 2016.

\bibitem[GPW23]{guo2023global}
Yan Guo, Benoit Pausader, and Klaus Widmayer.
\newblock Global axisymmetric euler flows with rotation.
\newblock {\em Inventiones mathematicae}, 231(1):169--262, 2023.

\bibitem[HLR{\etalchar{+}}19]{hickman2019higher}
Jonathan Hickman, Felipe Linares, Oscar~G Riano, Keith~M Rogers, and James
  Wright.
\newblock On a higher dimensional version of the benjamin--ono equation.
\newblock {\em SIAM Journal on Mathematical Analysis}, 51(6):4544--4569, 2019.

\bibitem[HXY18]{he2018global}
Ling-Bing He, Li~Xu, and Pin Yu.
\newblock On global dynamics of three dimensional magnetohydrodynamics:
  nonlinear stability of {A}lfv{\'e}n waves.
\newblock {\em Annals of PDE}, 4(1):1--105, 2018.

\bibitem[JW26]{JurjaWidmayer2026}
Catalina Jurja and Klaus Widmayer.
\newblock Long-time stability of a stably stratified rest state in the inviscid
  {2D} {Boussinesq} equation.
\newblock {\em Archive for Rational Mechanics and Analysis}, 250:26, 2026.

\bibitem[KLT14]{koh2014strichartz}
Youngwoo Koh, Sanghyuk Lee, and Ryo Takada.
\newblock Strichartz estimates for the euler equations in the rotational
  framework.
\newblock {\em Journal of Differential Equations}, 256(2):707--744, 2014.

\bibitem[Kno25a]{knobel2025suppression}
Niklas Knobel.
\newblock Suppression of fluid echoes and sobolev stability threshold for 2d
  dissipative fluid equations around couette flow.
\newblock {\em arXiv preprint arXiv:2505.23391}, 2025.

\bibitem[Kno25b]{knobel2025transition}
Niklas Knobel.
\newblock Transition and stability of 3d mhd around couette flow.
\newblock {\em arXiv preprint arXiv:2511.14609}, 2025.

\bibitem[Kno26]{Knobelideal}
Niklas Knobel.
\newblock Ideal {M}agnetohydrodynamics {A}round {C}ouette {F}low: {L}ong {T}ime
  {S}tability and {V}orticity-{C}urrent {I}nstability.
\newblock {\em Nonlinear Analysis}, 262:113937, 2026.

\bibitem[Ko24]{ko2024global}
Haram Ko.
\newblock Global axisymmetric solutions for navier-stokes equation with
  rotation uniformly in the inviscid limit.
\newblock {\em arXiv preprint arXiv:2409.17528}, 2024.

\bibitem[KZ23]{knobel2023echoes}
Niklas Knobel and Christian Zillinger.
\newblock On echoes in magnetohydrodynamics with magnetic dissipation.
\newblock {\em Journal of Differential Equations}, 367:625--688, 2023.

\bibitem[KZ25]{knobel2025sobolev}
Niklas Knobel and Christian Zillinger.
\newblock On the sobolev stability threshold for the 2d mhd equations with
  horizontal magnetic dissipation.
\newblock {\em Journal of Nonlinear Science}, 35(3):1--36, 2025.

\bibitem[Leh54]{lehnert1954magnetohydrodynamic}
B~Lehnert.
\newblock Magnetohydrodynamic waves under the action of the coriolis force.
\newblock {\em Astrophysical Journal, vol. 119, p. 647}, 119:647, 1954.

\bibitem[Leh55]{lehnert1955magnetohydrodynamic}
B~Lehnert.
\newblock Magnetohydrodynamic waves under the action of the coriolis force. ii.
\newblock {\em Astrophysical Journal, vol. 121, p. 481}, 121:481, 1955.

\bibitem[Lis20]{liss2020sobolev}
Kyle Liss.
\newblock On the {S}obolev stability threshold of 3{D} {C}ouette flow in a
  uniform magnetic field.
\newblock {\em Comm. Math. Phys.}, 377(2):859--908, 2020.

\bibitem[Mel24]{Melinand2024}
Benjamin Melinand.
\newblock Dispersive estimates for nonhomogeneous radial phases: An application
  to weakly dispersive equations and water wave models.
\newblock {\em Journal of Functional Analysis}, 286(1):110204, 2024.

\bibitem[NFV25a]{navarro2025nonlinear}
V{\'\i}ctor Navarro-Fern{\'a}ndez and David Villringer.
\newblock Nonlinear instability for the 3d mhd equations around the
  taylor-couette flow.
\newblock {\em arXiv preprint arXiv:2510.12707}, 2025.

\bibitem[NFV25b]{navarro2025spectral}
V{\'\i}ctor Navarro-Fern{\'a}ndez and David Villringer.
\newblock Spectral instability in the smooth {P}onomarenko dynamo.
\newblock {\em arXiv preprint arXiv:2509.19201}, 2025.

\bibitem[PW18]{pusateri2018global}
Fabio Pusateri and Klaus Widmayer.
\newblock On the global stability of a beta-plane equation.
\newblock {\em Analysis \& PDE}, 11(7):1587--1624, 2018.

\bibitem[RT24]{ren2024global}
Xiao Ren and Gang Tian.
\newblock Global solutions to the euler-coriolis system.
\newblock {\em arXiv preprint arXiv:2405.18390}, 2024.

\bibitem[RWZ21]{ren2021long}
Siqi Ren, Dongyi Wei, and Zhifei Zhang.
\newblock Long time behavior of {A}lfv{\'e}n waves in flowing plasma: the
  destruction of the magnetic island.
\newblock {\em SIAM Journal on Mathematical Analysis}, 53(5):5548--5579, 2021.

\bibitem[RZ17]{ren2017linear}
Siqi Ren and Weiren Zhao.
\newblock Linear damping of {A}lfv{\'e}n waves by phase mixing.
\newblock {\em SIAM Journal on Mathematical Analysis}, 49(3):2101--2137, 2017.

\bibitem[Tay63]{taylor1963magneto}
JB~Taylor.
\newblock The magneto-hydrodynamics of a rotating fluid and the earth’s
  dynamo problem.
\newblock {\em Proceedings of the Royal Society of London. Series A.
  Mathematical and Physical Sciences}, 274(1357):274--283, 1963.

\bibitem[WZ17]{wei2017global}
Dongyi Wei and Zhifei Zhang.
\newblock Global well-posedness of the {MHD} equations in a homogeneous
  magnetic field.
\newblock {\em Analysis \& PDE}, 10(6):1361--1406, 2017.

\bibitem[ZZ24]{zhao2024asymptotic}
Weiren Zhao and Ruizhao Zi.
\newblock Asymptotic stability of {C}ouette flow in a strong uniform magnetic
  field for the {E}uler-{MHD} system.
\newblock {\em Archive for Rational Mechanics and Analysis}, 248(3):47, 2024.

\bibitem[ZZZ21]{zhai2021long}
Cuili Zhai, Zhifei Zhang, and Weiren Zhao.
\newblock Long-time behavior of {A}lfv{\'e}n waves in a flowing plasma:
  generation of the magnetic island.
\newblock {\em Archive for Rational Mechanics and Analysis}, 242(3):1317--1394,
  2021.

\end{thebibliography}

\end{document}